\documentclass[pdflatex,sn-basic,Numbered]{sn-jnl}

\usepackage{graphicx}%
\usepackage{multirow}%
\usepackage{amsmath,amssymb,amsfonts}%
\usepackage{amsthm}%
\usepackage{mathrsfs}%
\usepackage[title]{appendix}%
\usepackage{xcolor}%
\usepackage{textcomp}%
\usepackage{manyfoot}%
\usepackage{booktabs}%
\usepackage{algorithm}%
\usepackage{algorithmicx}%
\usepackage{algpseudocode}%
\usepackage{listings}%

\usepackage{ulem}
\usepackage{graphicx}
\usepackage{amssymb,xfrac}
\usepackage{enumerate,color}
\usepackage{hyperref}
\usepackage{array}
\usepackage{psfrag}
\usepackage{pgf,tikz}
\usepackage{pgfplots}
\usetikzlibrary{arrows}
\usepackage{mathtools}
\pgfplotsset{compat=1.18}
\PassOptionsToPackage{hyphens}{url}
\usepackage{xurl}
\theoremstyle{thmstyleone}%
\newtheorem{lemma}{Lemma}

\theoremstyle{thmstyletwo}%

\theoremstyle{thmstylethree}%

\newcommand{\modulo}[1]{{\left|#1\right|}}
\newcommand{\norma}[1]{{\left\|#1\right\|}}

\newcommand{\R}{\mathbb{R}}

\newcommand{\N}{\mathbb{N}}

\newcommand{\Z}{\mathbb{Z}}

\renewcommand{\epsilon}{\varepsilon}
\renewcommand{\phi}{\varphi}
\renewcommand{\theta}{\vartheta}

\newcommand{\brho}{\boldsymbol{\rho}}

\newcommand{\Sim}{\mathcal{S}}
\newcommand{\F}{\mathbf{F}}
\newcommand{\dt}{{\Delta t}}
\newcommand{\dx}{{\Delta x}}

\newcommand{\be}{\begin{equation}}
\newcommand{\ee}{\end{equation}}

\definecolor{ffqqqq}{rgb}{1.,0.,0.}
\definecolor{uuuuuu}{rgb}{0.26666666666666666,0.26666666666666666,0.26666666666666666}

\begin{document}

\title[Discrete adjoint gradient computation for multiclass traffic flow models on road networks]{Discrete adjoint gradient computation for multiclass traffic flow models on road networks}

\author[1]{\fnm{Paola} \sur{Goatin}}
\email{paola.goatin@inria.fr}

\author[2]{\fnm{Axel} \sur{Klar}}
\email{klar@rptu.de}

\author*[1,2]{\fnm{Carmen} \sur{Mezquita-Nieto}}\email{carmen.mezquitanieto@math.rptu.de}

\affil[1]{\orgdiv{ACUMES}, \orgname{Universit\'e C\^ote d'Azur, Inria, CNRS, LJAD}, \orgaddress{\street{2004 route des Lucioles - BP 93}, \city{Sophia Antipolis Cedex}, \postcode{06902}, \country{France}}}

\affil[2]{\orgdiv{AGTM, Fachbereich Mathematik}, \orgname{RPTU Kaiserslautern}, \orgaddress{\street{Gottlieb-Daimler-Str. 48}, \city{Kaiserslautern}, \postcode{67663}, \country{Germany}}}

\abstract{This paper applies a discrete adjoint gradient computation method for a multi-class traffic flow model on road networks. Vehicle classes are characterized by their specific velocity functions, which depend on the total traffic density, resulting in a coupled hyperbolic system of conservation laws. The system is discretized using a Godunov-type finite volume scheme based on demand and supply functions, extended to handle complex junction coupling conditions --such as merges and diverges-- and boundary conditions with buffer lengths to account for congestion spillback. The optimization of different travel-related performance metrics, including total travel time and total travel distance, is formulated as a constrained minimization problem and is accomplished through the use of an adjoint gradient approach, allowing for an efficient computation of sensitivities with respect to the chosen time-dependent control variables. Numerical simulations on a sample network demonstrate the efficiency of the proposed framework, particularly as the number of control parameters increases. This approach provides a robust and computationally efficient solution, making it suitable for large-scale traffic network optimization.}


\maketitle

\section{Introduction}
\label{sec_art:introduction}

Due to the rapid increase in the demand for efficient control of traffic flow on road networks, there is a strong need for models that can accurately and efficiently capture the behavior of the dynamics and interactions between different types of vehicles, while also providing robust tools for optimization and control. 

Since the introduction of the Lighthill-Whitham-Richards (LWR) model~\cite{Light-Whit, Richards} and its extensions to consider several classes of vehicles on road networks, starting from~\cite{Benzoni-Colombo, WongWong2002}, there have been significant advances leading to more refined descriptions capturing the interactions between these classes, including for example creeping behaviour~\cite{FanWork2015, Gashaw2018}. Specifically, the extension to road networks incorporates the dynamics of junctions, yielding complex models that can accurately describe urban transportation systems~\cite{daganzo1995cell, garavello2010review, Joumaa2023}.

These advances naturally lead to the definitions of optimal control problems that can leverage these models to optimize traffic assignment. For an appropriate setting, it is important to specify the factors influencing traffic flow, establish an objective function that reflects the desired outcomes, and impose realistic constraints. Instead of relying on computationally expensive finite differences for the calculation of the gradient of the objective function, the adjoint method derives the gradient directly from the equations of the model; see, among others, \cite{GGK2016, Gugat2005, Reilly2015, SamaranayakeSO-DTA}. The choice of using the adjoint method over finite differences makes the difference in cases with a potential high number of control variables, as in the case in large road networks of time-dependent controls. There exist two main strategies for computing the adjoint gradient: one may derive the continuous adjoint equations from the underlying model, and then discretize it, or alternatively, one can formulate the adjoint equations directly on the discretized model, see~\cite{Gugat2005, HertyKirchnerKlar2007} for a comparative description of the two approaches. We remark that, since no convergence result is currently available for finite volume approximations in the multi-class setting, the equivalence of the above methods cannot be addressed in this case.
In this work, we therefore rely uniquely on the \textit{first discretize, then optimize} technique, relying on the (discrete) multi-class traffic flow model on networks proposed in~\cite{Joumaa2023, Joumaa2024}.
The adjoint gradient computation is then carefully implemented to take advantage of the sparsity of the involved lower triangular matrices~\cite{Reilly2015}, providing a computationally efficient mean of solving PDE-constrained traffic optimization problems by gradient descent methods.

The article is organized as follows: in Section~\ref{sec_art:multiclass} the multi-class traffic flow model and its discretization are introduced, specifying the dynamics at junction nodes and the inflow and outflow boundary conditions. In Section~\ref{sec_art:cost}, we describe the different cost functionals associated to road network performance, which we aim to optimize. Section~\ref{sec_art:DAG} details the gradient computation and structure of the adjoint system. In Section~\ref{sec_art:num_tests}, some numerical results are presented for a sample network with two vehicle classes to illustrate the effectiveness of the approach in optimizing the total travel time (TTT) and total travel distance (TTD) by controlling the distribution of each class of vehicles at road junctions. Finally, Section~\ref{sec_art:conclusions} discusses the results and future perspectives.

\section{A multi-class traffic flow model on networks with general speed function}
\label{sec_art:multiclass}

Multi-class macroscopic traffic models usually consist in systems of conservation equations, one for each vehicles class, which is characterized by its mean velocity function. In turn, class-specific velocities depend on all class densities, thus coupling the equations. We refer to~\cite{Logghe2008} for an overview of available descriptions.

In this paper, for the sake of simplicity but without restriction, we refer to the model class developed in~\cite{Benzoni-Colombo, FanWork2015, Joumaa2024, WongWong2002}, thus considering $N\times N$ systems of conservation laws 
\begin{equation}
    \brho_t + \textbf{F} (\brho)_x =0, \label{eq:CL}
\end{equation}
where 
\begin{equation}
    \brho = (\rho^1,\ldots,\rho^N)^T, \quad \F(\brho) = (v_1(r)\rho^1,\ldots,v_N(r)\rho^N)^T, \text{ and } r=\sum_{c=1}^N\rho^c.\label{eq:CLgeneral}
\end{equation}
Above, for $c=1,\ldots,N$, $\rho^c=\rho^c(x,t)$ denotes the vehicle class density and $v_c$ denotes the class-specific velocity function, which depends on the total traffic density. 
We assume
\begin{equation}
    v'_c(r)\leq 0 \label{eq:hypothesis1}
\end{equation}
and
\begin{equation}
    v_c(0) = V_c,\quad \exists\, R_c \text{ s.t. } v_c(r) = 0 \text{ for } r\geq R_c \label{eq:hypothesis2}
\end{equation}
for $c=1,\ldots,N$, where $R_c$ is the maximal density for each class. 
Without loss of generality, we also assume that
\begin{equation}
    V_1\geq \ldots \geq V_N. \label{eq:speeds}
\end{equation}
System \eqref{eq:CL} is then defined on the $n$-dimensional polyhedron
\begin{equation*}
    \Sim = \left\{ \brho\in\R^N \colon 0\leq \rho^c \leq R_c \text{ and }\sum_{c=1}^N\rho^c \leq R \right\},
\end{equation*}
where $R\coloneqq\max\{R_1,\ldots,R_N\}$. 
As shown in~\cite[Theorem 2.1]{Joumaa2024}, the system is hyperbolic. 

Following~\cite{BrianiCristiani2014, Burger2008, Joumaa2023}, we
approximate system~\eqref{eq:CL} by the Godunov-type finite volume scheme~\cite{Godunov}
\begin{equation} 
    \rho_j^{c,\nu+1} = \rho_j^{c,\nu} - \dfrac{\dt}{\dx} \left[ F^{c,\nu}_{j+1/2} - F^{c,\nu}_{j-1/2}\right], \quad j\in\Z,\,\nu\in\N,\,c=1,\ldots,N, \label{eq:FV}
\end{equation}
with
\begin{equation}  
    F^{c,\nu}_{j+1/2} \coloneqq \frac{\rho_j^{c,\nu}}{r_j^\nu}\min\left\{ D^c(r_j^\nu), S^c(r_{j+1}^\nu)\right\}, \label{eq:MPgodunov}
\end{equation}
where $D^c$ and $S^c$ are respectively the demand and supply functions~\cite{Lebacque1996GodunovScheme} of the total density relatively to the $c$-th class speed defined by setting $Q^c(r):=v_c(r)r$, $r^c_{cr}:=\arg\max_{r} Q^c (r)$ and 
\begin{align}
     D^c(r) &\coloneqq Q^c(\min\{r,r^c_{cr}\}), \label{eq:MPdemand} \\
    S^c(r) &\coloneqq Q^c(\max\{r,r^c_{cr}\}). \label{eq:MPsupply}
\end{align}
See also~\cite[Eq. (4)]{LevinBoyles2016}.
Other formulations of the demand and supply functions are presented in the literature, see~\cite{FanWork2015, Herty2006}, but they may lead to non-physical numerical oscillations. 

Under a suitable Courant-Friedrichs-Lewy (CFL)~\cite{CFL1928} stability condition, the scheme~\eqref{eq:FV} - \eqref{eq:MPgodunov} guarantees the invariance of the set $\Sim$ and therefore the fulfillment of the natural density bounds~\cite[Proposition 2.4]{Joumaa2024}:
\begin{lemma} 
    \label{lem:Gboundgeneral}
    Under the CFL condition
    \begin{equation}
        \max\left\{\max_c{\norma{v_c}}_\infty, \max_c{\norma{Q_c'}_\infty}\right\} \Delta t \leq \Delta x, \label{eq:cflUWgeneral}
    \end{equation}
    for any initial data $\brho_0\in\Sim$, the approximate solutions computed by the Godunov scheme \eqref{eq:FV} - \eqref{eq:MPgodunov} satisfy the following uniform bounds:
    \begin{equation*}
        \brho_j^\nu\in\Sim \quad \forall j\in\Z,\,\nu\in\N.
    \end{equation*}
\end{lemma}

\subsection{Coupling conditions at junctions} 
\label{sec_art:multiclass:junctions}

In the following, we provide a constructive definition of the solution of the Riemann Problem for general junctions (see \cite{coclite2005traffic, DelleMonachePriorityRS, GaravelloPiccoli2006}), i.e. for the case with $n$ incoming roads and $m$ outgoing roads for $N$ classes of vehicles, as detailed in~\cite{Joumaa2024}. The junction must fulfill the property of mass-conservation at junctions, i.e. the sum of the fluxes $\hat{\gamma}_{i,out}^{c,\nu} $ entering a junction has to be equal to the sum of the fluxes $\hat{\gamma}_{j,in}^{c,\nu}$ leaving it, for any fixed class $c$ of vehicles and iteration $\nu$,
\begin{align}
    \sum_{i=1}^n \hat{\gamma}_{i,out}^{c,\nu} &= \sum_{j=n+1}^{n+m} \hat{\gamma}_{j,in}^{c,\nu}. \label{eq:mass_conserv}
\end{align}
For this Riemann Problem, we need to consider the subsequent rules to be able to determine a unique solution on the network:
\begin{itemize}
    \item[] \textbf{(A)} There is a class-specific distribution matrix $A^c$, $c=1,\ldots,N$,
    \begin{equation*}
        A^c =
        \begin{pmatrix}
            \alpha_{n+1,1}^c & \cdots & \alpha_{n+1,n}^c\\
            \vdots & \ddots & \vdots\\
            \alpha_{n+m,1}^c & \cdots & \alpha_{n+m,n}^c
        \end{pmatrix},
    \end{equation*}
    where the coefficients $\alpha_{ji}^c\geq 0$ represent the percentage of vehicles of class $c$ going from road $i\in\{1,\ldots,n\}$ to road $j\in\{n+1,\ldots,n+m\}$. In particular, to ensure mass conservation, we require that for any $i\in\{1,\ldots,n\}$, it holds
    \begin{equation*}
        \sum_{j=n+1}^{n+m} \alpha_{ji}^c = 1.
    \end{equation*}

    \item[] \textbf{(B)} Respecting \textbf{(A)}, the sum of the fluxes through the junction is maximized, i.e.
    \begin{equation}
        \sum_{i=1}^n \hat{\gamma}^{c,\nu}_{i,out}  \quad \text{is maximal,}
    \end{equation}
    and
    \begin{equation}
        \hat{\gamma}_{j,in}^{c,\nu} = \sum_{i=1}^n \alpha_{ji}^c \hat{\gamma}_{i,out}^{c,\nu}, \quad j=n+1,\ldots,n+m.
    \end{equation}
\end{itemize}
Further assumptions are often required to ensure uniqueness of the solution, such as specific conditions on the distribution coefficients~\cite{GaravelloPiccoli2006} or the introduction of priority parameters for the incoming roads~\cite{DelleMonachePriorityRS}.
\begin{figure}[!t] % Figure 1
	\centering
	\includegraphics[]{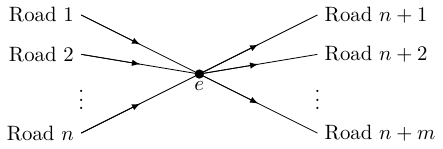} % general junction
    \caption{Structure of a general $n\times m$ junction $e$,}
    \label{fig:junctions:nxm:general_junction}
\end{figure}

In the case of simple junctions, i.e. for $1 \times 1$, $M \times 1$ merges and $1 \times M$ diverges, 
we can write explicitly  the solutions of the Riemann Problem  for $N$ classes of vehicles. Below, the demand and supply function are given by equations~\eqref{eq:MPdemand} and~\eqref{eq:MPsupply}, but they may refer to different fluxes $Q^c=Q^c_\ell$ on each road, thus depending on the index $\ell=1,\ldots,M+1$.

\subsubsection{Solution for \texorpdfstring{$1\times 1$}{1x1} junctions} 
\label{sec_art:multiclass:junctions:1x1}

The original numerical flux \eqref{eq:MPgodunov} is applied to the $1 \times 1$ junction between an incoming road $\ell=1$ and an outgoing road $\ell=2$, with corresponding discretized densities $\rho_{\ell,j}^{c,\nu}$, as
\begin{equation*}
    \hat{\gamma}^{c,\nu}_{1,out}=\hat{\gamma}^{c,\nu}_{2,in}=\frac{\rho_{1,N_1}^{c,\nu}}{r_{1,N_1}^{\nu}} \min\left\{D^c_1(r_{1,N_1}^{\nu}),S^c_2(r_{2,1}^{\nu})\right\}, \qquad c=1,\ldots,N,
\end{equation*}
where $j=N_1$ denotes the last cell of the incoming road and $j=1$ the first cell of the outgoing road.

\subsubsection{Solution for \texorpdfstring{$M\times 1$}{Mx1} (merge) junctions} \label{sec_art:multiclass:junctions:Mx1}

In this case, we need to 
introduce a class-specific priority vector $P^c=(p_1^c,\ldots,p_M^c)\in\R^M$, $p_i^c\geq 0$, $\sum_{i=1}^M p_i^c=1$, so that $\hat{\gamma}^{c,\nu}_{i,out}=p_i^c \,\hat{\gamma}^{c,\nu}_{M+1,in}$ for $i=1,\ldots,M$.

The corresponding fluxes are:
\begin{align} \label{eq:merge_junction}
    \begin{aligned}
        \begin{split}
            \hat{\gamma}^{c,\nu}_{i,out} &=\dfrac{\rho_{i,N_i}^{c,\nu}}{r_{i,N_i}^{\nu}} \min \Big\{ D^c_i(r_{i,N_i}^{\nu}), \max \big\{p_i^c S^c_{M+1}(r_{M+1,1}^{\nu}),\\
            &\qquad\qquad\qquad\qquad\qquad\qquad\quad\; S^c_{M+1}(r_{M+1,1}^{\nu})- \sum_{j\not= i} D^c_j(r_{j,N_j}^{\nu}) \big\} \Big\},
        \end{split} \\
        \hat{\gamma}^{c,\nu}_{M+1,in} &=\sum_{i=1}^M\hat{\gamma}^{c,\nu}_{i,out},
    \end{aligned}
\end{align}
for $c=1,\ldots,N$ and $i=1,\ldots,M$. 
\begin{figure}[!ht] % Figure 2
	\centering
	\includegraphics[]{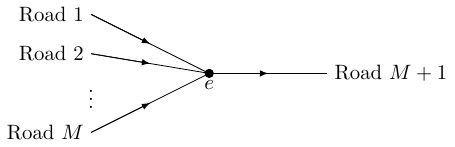} % merging junction
    \caption{Structure of a $M\times 1$ junction $e$.}
    \label{fig:junctions:Mx1:merging_junction}
\end{figure}

\subsubsection{Solution for \texorpdfstring{$1\times M$}{1xM} (diverge) junctions} \label{sec_art:multiclass:junctions:1xM}

In this case, the class-specific distribution matrix $A^c$ takes the form
\begin{equation*}
    A^c = (\alpha_2^c,\ldots,\alpha_{M+1}^c)^T,
\end{equation*}
with $\alpha_i^c \geq 0$, $\sum_{i=2}^{M+1} \alpha_i^c=1$. For the FIFO (first-in-first-out) version~\cite{GaravelloPiccoli2006}, the fluxes at the junction are defined as
\begin{align} \label{eq:FIFO_div_junction} 
    \begin{aligned}
        \hat{\gamma}^{c,\nu}_{i,in} &= \alpha_i^c \hat{\gamma}^{c,\nu}_{1,out}, \\
        \hat{\gamma}^{c,\nu}_{1,out} &= \frac{\rho_{1,N_1}^{c,\nu}}{r_{1,N_1}^{\nu}} \min \left\{ D_1^c(r_{1,N_1}^{\nu}),\frac{S_2^c(r_{2,1}^{\nu})}{\alpha_2^c},\ldots,\frac{S_{M+1}^c(r_{M+1,1}^{\nu})}{\alpha_{M+1}^c} \right\},
    \end{aligned}
\end{align}
for $c=1,\ldots,N$ and $i=2,\ldots,M+1$.
The limitation of this description consists in the fact that the whole traffic is blocked as soon as one outgoing road is fully congested. Since this may not be realistic in some situations, we also consider   
the non-FIFO version~\cite{GGK2016, HertyKlar2003introductionnonFIFO, LebacqueKhoshyaran2002} below:
\begin{align} \label{eq:nonFIFO_div_junction}
    \begin{aligned}
        \hat{\gamma}^{c,\nu}_{i,in} &= \dfrac{\rho_{1,N_1}^{c,\nu}}{r_{1,N_1}^{\nu}} \min \left\{\alpha_i^c D^c_1(r_{1,N_1}^{\nu}), S^c_i(r_{i,1}^{\nu}) \right\}, \\[5pt]
        \hat{\gamma}^{c,\nu}_{1,out} &= \sum_{i=2}^{M+1}\hat{\gamma}^{c,\nu}_{i,in},
    \end{aligned}
\end{align}
for $c=1,\ldots,N$ and $i=2,\ldots,M+1$. For the numerical experiments presented in Section \ref{sec_art:num_tests}, the diverging junction used will be a FIFO \eqref{eq:FIFO_div_junction} type.
\begin{figure}[!ht] % Figure 3
	\centering
	\includegraphics[]{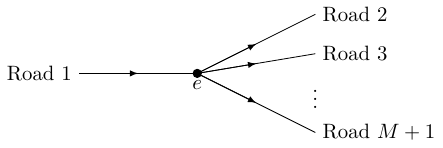} % diverging junction
    \caption{Structure of a $1\times M$ junction $e$.}
    \label{fig:junctions:1xM:diverging_junction}
\end{figure}

\subsection{Boundary conditions}
\label{sec_art:multiclass:BC}

Let $F_{in}=(F_{in}^1,\ldots,F_{in}^N)^T$ and $F_{out}=(F_{out}^1,\ldots,F_{out}^N)^T$, with $F_{in}^c,F_{out}^c\geq 0$ for all $c=1,\ldots,N$, be the inflow boundary conditions at an origin node and the outflow boundary condition at a destination node, respectively.

In order not to miss vehicle counting when congestion spills back to the incoming nodes of the network, a buffer of length $l_c(t)$ at the beginning of the entrance road is considered (see e.g.~\cite{GGK2016}). The length of the buffer must be updated at every time step in order to choose the correct demand function to calculate the flux at the beginning of the buffer.

From \eqref{eq:MPgodunov}, \eqref{eq:merge_junction}, we can define
\begin{align*}
    \hat{\gamma}^{c,\nu}_{\ell,in} & = \min \bigg\{D^c_{\ell}(l_{\ell}^{c,\nu}), \max \Big\{ \dfrac{1}{N} S^c_{\ell}(r_{\ell,1}^{\nu}),S^c_{\ell}(r_{\ell,1}^{\nu})- \sum_{g\not=c}D^g_{\ell}(l^{g,\nu}_{\ell}) \Big\} \bigg\}, \\
    \hat{\gamma}^{c,\nu}_{\ell,out} &= \min \left\{  \dfrac{\rho_{\ell,N_{\ell}}^{c,\nu}}{r_{\ell,N_{\ell}}^{\nu}}D^c_{\ell}(r_{\ell,N_{\ell}}^{\nu}),F_{\ell,out}^c\right\},
\end{align*}
and 
\begin{equation*}
    D^c_{\ell}(l^{c,\nu}_{\ell})=
    \begin{cases}
        Q^c_{\ell}(r^c_{\ell,cr}) & \text{ if } l^{c,\nu}_{\ell}>0,\\
        F_{\ell,in}^c & \text{ if } l^{c,\nu}_{\ell}=0.
    \end{cases}
\end{equation*}
To update the length at every time step, we set
\begin{equation*}
    l_{\ell}^{c,\nu+1}=\max \{l^{c,\nu}_{\ell}+\Delta t (F_{\ell,in}^c-\hat{\gamma}^{c,\nu}_{\ell,in}), 0 \},
\end{equation*}

\section{Cost functionals} 
\label{sec_art:cost}

Cost functionals assess and regulate the performance of the traffic network searching to achieve different goals, which can be reducing time delays, improving traffic flow, etc. In the multi-class case, this approach evaluates as well the interactions and behaviors of different vehicle classes. To study these interactions and how to better regulate them via suitable control actions, we consider two metrics, the total travel time and the total travel distance, which will be later used in the optimization procedures analyzed in Section~\ref{sec_art:num_tests}. It is possible to contemplate other metrics, such as the network throughput, i.e. the total number of vehicles of that class leaving the network~\cite{Joumaa2023}, the $\text{CO}_2$ emissions~\cite{Joumaa2024} or the mean speed~\cite{PiacentiniCicic2019}.

The total travel time (TTT) of a given class $c$, i.e. the space and time integral of the corresponding density, is defined as
\begin{equation*}
    \text{TTT}(\rho^c) = \sum_{\ell} \int_0^{T_f} \int_{a_{\ell}}^{b_{\ell}} \rho_{\ell}^c(x,t) dx dt + \sum_i \int_0^{T_f} l_i^c(t) dt,
\end{equation*}
where $[a_{\ell},b_{\ell}]\subset\R$ is the space interval corresponding to road $\ell$, $[0,T_f]$ the time interval considered and $i$ the origin nodes index. The discrete counterpart is then defined as
\begin{equation}
    \text{TTT}(\rho^c) = \Delta t \Delta x \sum_{\nu=0}^{T_f} \sum_{\ell} \sum_{j=0}^{N_\ell} \rho^{c,\nu}_{\ell,j} + \Delta t \sum_{i}\sum_{\nu=0}^{T_f} l_i^{c,\nu},
    \label{eq:TTT}
\end{equation}
where $N_\ell$ denotes the number of cells in the road segment $\ell$. \\
The total travel distance (TTD) of the vehicle class $c$ corresponds to the total distance travelled by the vehicles,
\begin{equation*}
    \text{TTD}(\rho^c) = \sum_{\ell} \int_0^{T_f} \int_{a_{\ell}}^{b_{\ell}} F_{\ell}^c(x,t) dx dt, \qquad F^c_{\ell}(x,t)\coloneqq \rho_{\ell}^c(x,t) v_{\ell,c}(x,t),
\end{equation*}
or, discretized,
\begin{equation}
    \text{TTD}(\rho^c) = \Delta t\Delta x\sum_{\nu=0}^{T_f} \sum_{\ell} \sum_{j=0}^{N_{\ell}} \rho_{\ell,j}^{c,\nu} v_{\ell,c}(r_{\ell,j}^{\nu}).
    \label{eq:TTD}
\end{equation}
By linearity, the global functionals are given by the sum of the corresponding class-specific quantities:
\begin{equation*}
    \text{TTT}(r) = \sum_{c=1}^N \text{TTT}(\rho^c), \qquad \text{TTD}(r) = \sum_{c=1}^N \text{TTD}(\rho^c).
\end{equation*}
The TTT and TTD are closely related. In uncongested settings with a constant speed throughout the network, the buffer $l_{\ell}^c$ vanishes and therefore \eqref{eq:TTT} and \eqref{eq:TTD} are linearly dependent, with the quotient $\text{TTD}/\text{TTT}$ providing a definition for the mean speed. However, in congested settings, the relationship between TTT and TTD becomes nonlinear, with smaller vehicle velocities leading to higher travel time values for the same travelled distance.

Minimizing one of these quantities typically results in an increase of the other: reducing delay and congestion to minimize the TTT may lead to vehicles having to take a longer route, therefore increasing the TTD, while shortening travel distances may lead to a higher TTT due to a denser traffic. To solve this drawback and find a balanced situation for both functionals, it is typical to consider multi-objective optimization problems, which will be treated in Section \ref{sec_art:num_tests:Pareto_front}.

\section{Discrete adjoint gradient computation} 
\label{sec_art:DAG}

Let $J(u,y)$ be a function to be minimized, where $y$ is the vector of state variables containing the approximations of the densities $\rho$ and the buffer lengths $l$, i.e., $y=(l,\rho)$, and $u$ is a vector containing all the control variables. Both $y$ and $u$ satisfy the underlying model equations, denoted by $E(u,y)=0$. 
In the specific case of the traffic flow on a road network, the control variables generally consist of a vector $u=(V,A,P)$, with
\begin{itemize}
    \item $V=\{V_{\ell}\}_{\ell}$, $V_{\ell} = (V_{\ell,1},\ldots,V_{\ell,N})$, where $V_{\ell,c}\in [0, V_c^{\max}]$ corresponds to the (possibly time-dependent) class-specific speed limit on road segment $\ell$ (these controls are used to implement class-specific variable speed limits on the roads of the network);

    \item $A=\{A_e\}_e$, $A_e=(A_e^{1},\ldots,A_e^{N})$, where $A_e^{c}$ is the (time-dependent) class-specific distribution matrix at junction $e$ (providing a routing control action specific for each class of vehicles);
    
    \item $P=\{P_e\}_e$, $P_e=(P_e^{1},\ldots,P_e^{N})^T$, where $P_e^{c}$ is the (time-dependent) class-specific priority vector on junction $e$ (which can be used to implement ramp metering or other regulations through traffic light signals). 
\end{itemize}
In this setting, a general traffic management problem can be rewritten as a PDE-constrained optimization problem:
\begin{equation}
    \min_{u\in\mathcal{U}} J(u,y), \qquad \text{subject to } E(u,y) = 0,\label{eq:OptProbl}
\end{equation}
where $\mathcal{U}$ denotes the set of admissible controls, to be specified depending on the situation.

\begin{figure}[!hb] % Figure 5
    \centering
    \includegraphics[width=0.49\textwidth]{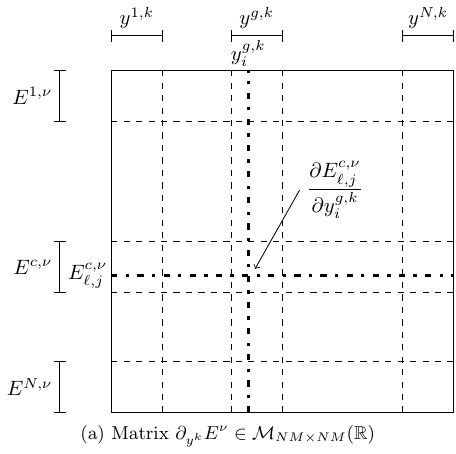}\hfill
    \includegraphics[width=0.48\textwidth]{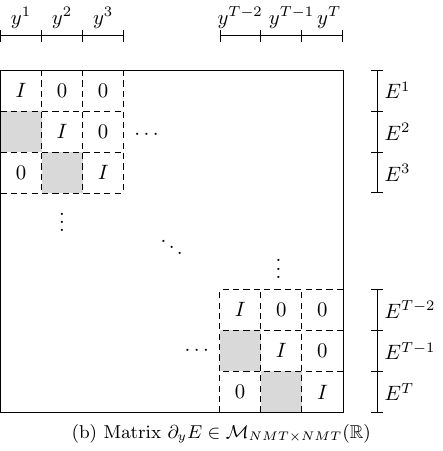}
    \caption{Partial derivative terms in $\partial_y E$ are ordered first by time and then by cell index. Blocks where time iteration $k\notin\{\nu-1,\nu\}$ are 0.} 
    \label{fig:DAG:struct_matrix_E_y}
\end{figure}

Assuming $E$ has a unique solution for each $u$, the reduced model $J(u,y(u))$ can be considered. The cost gradient is then computed as
\begin{equation*}
    \frac{d}{du}J(u,y(u)) = \partial_u J + \partial_y J\, \frac{dy}{du}\,.
\end{equation*}
It can be deduced from the equation $E(u,y(u))=0$ that
\begin{equation*}
    \frac{d}{du}E(u,y(u)) = \partial_u E + \partial_y E\, \frac{dy}{du} = 0 \Rightarrow \frac{dy}{du} = -(\partial_y E)^{-1}\,\partial_u E.
\end{equation*}
Therefore, substituting in the cost gradient,
\begin{align}
    \frac{d}{du}J(u,y(u)) &= \partial_u J - \partial_y J\,(\partial_y E)^{-1}\,\partial_u E \label{eq:costgradient}\\
    &= \partial_u J + \xi^T \partial_u E, \nonumber
\end{align}
where $\xi$, the adjoint variable, is the solution to the adjoint equation
\begin{equation}
    (\partial_y E)^T\xi = -(\partial_y J)^T. \label{eq:adjointeq}
\end{equation}
To compute the gradient, the following partial derivatives must then be obtained:
\begin{equation}
  \frac{\partial J}{\partial y}, \quad \frac{\partial J}{\partial u}, \quad \frac{\partial E}{\partial y}, \quad \frac{\partial E}{\partial u}. \label{eq:derivatives_costgradient}
\end{equation}
Let now $\nu=1,\ldots,T$, where $T=T_f$ is the number of time steps. We denote by $\mathcal{L}$ the set of indexes of all links, and by
$\mathcal{O}\subset \mathcal{L}$ the subset of indexes of incoming links with buffers. Then, the state variables for each road and class are
\begin{itemize}
    \item $(l_{\ell}^{c,\nu},\rho_{\ell,1}^{c,\nu},\ldots,\rho_{\ell,N_\ell}^{c,\nu})$ for $\ell\in\mathcal{O}$,
    \item $(\rho_{\ell,1}^{c,\nu},\ldots,\rho_{\ell,N_\ell}^{c,\nu})$ for $\ell\in\mathcal{L}\setminus\mathcal{O}$.
\end{itemize}
Let $M:=\sum_{\ell\in\mathcal{L}}N_\ell + \modulo{\mathcal{O}}$ be the total number of cells in the network. Additionally, define $U$ as the total number of control points considered, with $U\coloneqq |\mathcal{U}|$. For each state variable $y$, the corresponding cost function $J$ is defined as
\begin{align*}
    J\colon \mathbb{R}^{NUT} \times \mathbb{R}^{NMT} &\to \mathbb{R}\\
    (u,y) &\mapsto J(u,y),
\end{align*}
where $y=(y^1,\ldots,y^T)$, $y^{\nu}=(y^{1,\nu},\ldots,y^{N,\nu})$. 
\begin{figure}[!t] % Figure 4
	\centering
	\includegraphics[]{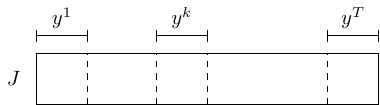} % vector Jy
    \caption{Structure of the row vector $\partial_y J\in\R^{NMT}$.}
    \label{fig:DAG:struct_vector_J_y}
\end{figure}
The update equations $E$ are defined as
\begin{align*}
    E\colon \mathbb{R}^{NUT} \times \mathbb{R}^{NMT} &\to \mathbb{R}^{NMT}\\
    (u,y) &\mapsto E(u,y) = 0,
\end{align*}
where each element of $E$ is defined as $E_{\ell}^{c,\nu}$, with
\begin{equation*}
    E_{\ell}^{c,\nu} = (E_{\ell,j}^{c,\nu})_{j=1,\ldots,N_{\ell}},\quad E^{\nu} = \left((E_{\ell}^{1,\nu})_{\ell\in\mathcal{L}},\ldots,(E_{\ell}^{N,\nu})_{\ell\in\mathcal{L}}\right),\quad E = (E^1,\ldots,E^T).
\end{equation*}
In particular, for $c=1,\ldots,N$, $\nu=1,\ldots,T$ and $\ell\in\mathcal{L}$,
\begin{equation}
    0 = E_{\ell,j}^{c,\nu}(u,y) \coloneqq
    \begin{cases}
        \rho_{\ell,1}^{c,\nu} - \rho_{\ell,1}^{c,\nu-1} + \lambda\, \left[F_{\ell,3/2}^{c,\nu-1} - \hat{\gamma}^{c,\nu-1}_{\ell,in}\right], & j=1,\\[5pt]
        \rho_{\ell,j}^{c,\nu} - \rho_{\ell,j}^{c,\nu-1} + \lambda \left[ F_{\ell,j+1/2}^{c,\nu-1} - F_{\ell,j-1/2}^{c,\nu-1} \right], & j=2,\ldots,N_{\ell}-1, \\[5pt]
        \rho_{\ell,N_{\ell}}^{c,\nu} - \rho_{\ell,N_{\ell}}^{c,\nu-1} + \lambda\left[\hat{\gamma}^{c,\nu-1}_{\ell,out} - F_{\ell,N_{\ell}-1/2}^{c,\nu-1}\right], & j=N_{\ell},
    \end{cases} 
    \label{eq:DiscModelEq_rho}
\end{equation}
and, for $\ell\in\mathcal{O}$,
\begin{equation}
    0 = E_{\ell,0}^{c,\nu}(u,y) \coloneqq l_{\ell}^{c,\nu} - \max \left\lbrace l_{\ell}^{c,\nu-1} + \Delta t (F_{\ell,in}^c-\hat{\gamma}^{c,\nu-1}_{\ell,in}), 0 \right\rbrace,
    \label{eq:DiscModelEq_l}
\end{equation}
with $\lambda=\Delta t/\Delta x$. The next step is to compute the partial derivatives in \eqref{eq:derivatives_costgradient}.
The definition of $\partial_y E$ shows a structure similar to the one proposed in~\cite{Reilly2015}. $E_{\ell,j}^{c,\nu}$ depends only on the previous time step, i.e. on the class density terms $\rho_{\ell,j-1}^{c,\nu-1}$, $\rho_{\ell,j}^{c,\nu-1}$, $\rho_{\ell,j+1}^{c,\nu-1}$ and the total density terms $r_{\ell,j-1}^{\nu-1}$, $r_{\ell,j}^{\nu-1}$, $r_{\ell,j+1}^{\nu-1}$ (as well as $l_{\ell}^{c,\nu-1}$ if $\ell\in\mathcal{O}$), and the current time step, with $\rho_{\ell,j}^{c,\nu}$ (and $l_{\ell}^{c,\nu}$ if $\ell\in\mathcal{O}$). In particular, the terms of $\partial_y E$ on the main diagonal are equal to 1, while those values over it are equal to 0.
Hence, the matrix $\partial_y E$ has a lower-triangular sparse structure and it is invertible, as seen in Figure \ref{fig:DAG:struct_matrix_E_y}, which is later used to compute the adjoint variable $\xi$, and its high sparsity reduces significantly the computational cost of the gradient. Due to this structure, solving the adjoint system \eqref{eq:adjointeq} can be done efficiently through a backward substitution of the upper triangular matrix $(\partial_y E)^T$. 

Obtaining the cost gradient \eqref{eq:costgradient} through the calculation of the adjoint variable $\xi$ requires at most $O((NMT)^2 + (NMT)(NUT))$ steps: finding the solution of an upper-triangular $NMT\times NMT$ system and a multiplication of two matrices of sizes $NMT\times NMT$ and $NMT\times NUT$. As long as $NUT>1$, this would be more efficient than calculating the cost gradient \eqref{eq:costgradient} without solving the adjoint system, i.e. determining it instead through a direct forward substitution, which in turn requires solving lower-triangular $NMT\times NMT$ systems $NUT$ times, with complexity $O((NMT)^2(NUT))$. Additional explanations regarding the computational costs for the one-class method can be found in \cite{Reilly2015}, including a reduction of the complexity for the adjoint gradient.

For a more detailed computation of the different derivatives of the discrete adjoint gradient, see Appendix \ref{app:DAG}.

\section{Numerical tests}
\label{sec_art:num_tests}

\begin{figure}[!hb] % Figure 6
	\centering
	\includegraphics[]{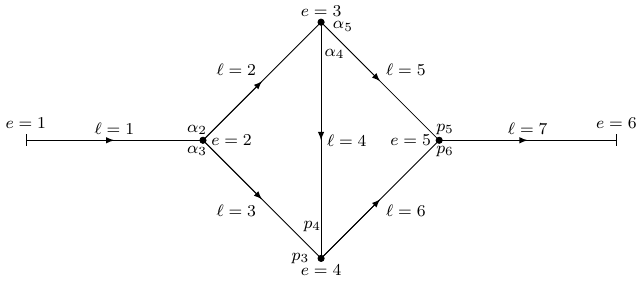} % test network
    \caption[]{Diagram of the test network showing the different roads and their associated parameters. Junctions $e=2,3$ diverge and junctions $e=4,5$ merge.}
    \label{fig:num_tests:example:network}
\end{figure}

A series of numerical tests are conducted for two  classes of vehicles moving along the sample network proposed in Figure \ref{fig:num_tests:example:network}. The goal is to observe how varying the possible control variables influence the metrics defined in Section \ref{sec_art:cost}. In particular,  two tests will be presented in this section. The first one compares the optimal minimal solutions of the TTT \eqref{eq:TTT} obtained through central finite differences and through an adjoint gradient computation, while the second test generates Pareto fronts for a multi-objective optimization of  TTT \eqref{eq:TTT} and TTD \eqref{eq:TTD}. In both tests, we will introduce time-dependent control variables to address the changing conditions of traffic flow over time. All simulations are run on a laptop with a $13^{\text{th}}$ generation Intel$^{\text{(R)}}$ Core$^{\text{(TM)}}$ i5-1335U 1.30 GHz processor.

For two classes of vehicles, the total density $r$ equals the sum of the densities of the two classes $\rho^1$ and $\rho^2$, i.e. $r=\rho^1+\rho^2$, with
\begin{equation*}
    F_{\ell}^c(\rho^1,\rho^2) = v_{\ell,c}(r)\rho^c \quad \text{and} \quad Q_{\ell}^c(r) = v_{\ell,c}(r)r, \quad c=1,2.
\end{equation*}
For the following results, we consider the Greenshields' velocity function~\cite{Greenshields1935} given by
\begin{equation*}
    v_{\ell,c} = V_{\ell,c}\left(1-\dfrac{r}{R_{\ell,c}}\right), \quad \ell=1,\ldots,7,
\end{equation*}
where $V_{\ell,1}\geq V_{\ell,2}$ for each $\ell$, assuming the first class travels generally faster than the second. This means that the critical density associated to a road $\ell$ and a class $c$ is $r_{\ell,cr}^c=R_{\ell,c}/2$, and thus $Q_{\ell,cr}^c\coloneqq Q_{\ell}^c(r_{\ell,cr}^c)=V_{\ell,c} R_{\ell,c}/4$. The lengths $L_{\ell}$, $\ell=1,\ldots,7$, of the roads are set to
\begin{equation*}
    L_1=L_2=L_4=L_6=L_7=3, \qquad L_3=L_5=6,
\end{equation*}
so that each possible route from the origin node $e=1$ to the destination node $e=6$ (routes 1-2-5-7, 1-2-4-6-7, and 1-3-6-7) has length 15. The maximum velocities and maximum total densities on each road $\ell$ and class $c$ are given by 
\begin{equation*}
    V_1 = (80,80,80,80,80,80,80), \qquad V_2 = (80,40,40,20,20,40,40),
\end{equation*}
and $R_{\ell,1} = R_{\ell,2} = R = 150$, respectively. For the definition of the boundary conditions, the inflow $F_{in}\coloneqq F_{1,in} = (F_{1,in}^1,F_{1,in}^2)$ is set as
\begin{equation*}
    F_{in} =
    \begin{cases}
        (Q_{1,cr}^1, \frac{2}{3}\,Q_{1,cr}^2) & \text{ if } t\leq T_f/2,\\
        (0,0) & \text{ if } t>T_f/2,
    \end{cases}
\end{equation*}
where $T_f=1$ denotes the simulation time horizon and $F_{out}\coloneqq F_{7,out} = (F_{7,out}^1,F_{7,out}^2)$ maximizes the outflow of the network, i.e.
\begin{equation*}
    F_{out} = (Q_{7,cr}^1, Q_{7,cr}^2).
\end{equation*}
In the numerical tests, we choose the splitting parameters $\alpha^c$ and $\beta^c$ as control variables in the numerical tests, since they directly determine the distribution of the flow of each vehicle class at the diverging junctions $e=2,3$ along the alternative routes of the network. From a practical perspective, these parameters can be interpreted as routing controls corresponding to real-life traffic management actions, where they can be used to divert vehicles away from potentially congested roads and thus reduce the formation of traffic jams. Therefore, the diverging parameters are suitable and effective for improving the overall network performance.

The priority parameters for the merging junctions $e=4,5$ are kept fixed for both classes $c=1,2$, with values
\begin{equation*}
    p_3^c=p=1/2, \quad p_5^c=q=1/3 \quad (\text{and } p_4^c=1-p=1/2, \quad p_6^c=1-q=2/3).
\end{equation*}
The initial conditions of the density and buffer length $l\coloneqq l_1=(l_1^1,l_1^2)$ on the roads are set equal to 0, and the number of gridpoints used to discretize each road is $N_{\ell}=10L_{\ell}$, $\ell=1,\ldots,7$. The Godunov scheme \eqref{eq:FV}-\eqref{eq:MPgodunov} is implemented everywhere, except on the first and last cell of every road, where the fluxes $\hat{\gamma}_{\ell,in}^{c,\nu}$ and $\hat{\gamma}_{\ell,out}^{c,\nu}$ are defined as in Section \ref{sec_art:multiclass:junctions}. In particular, junctions $e=2,3$ follow the FIFO fluxes \eqref{eq:FIFO_div_junction}, while \eqref{eq:merge_junction} is applied for junctions $e=4,5$.

Let us concatenate the vectors $(\rho_{\ell,i}^c)_i$, $c=1,2$, and $(r_{\ell,i})_i$, $\ell=1,\ldots,7$, $i=1,\ldots,N_{\ell}$ into vectors of size $\sum_{i=1}^7 N_{\ell}$ defined as 
\begin{align*}
    \rho^c &= (\rho_1^c,\ldots,\rho^c_{N_1+\cdots+N_7}) = ((\rho^c_{1,i})_i,\ldots,(\rho^c_{7,i})_i),\\
    r &= (r_1,\ldots,r_{N_1+\cdots+N_7}) = ((r_{1,i})_i,\ldots,(r_{7,i})_i),
\end{align*}
where $\rho_j^c=\rho_{\ell,i}^c$, $r_j=r_{\ell,i}$, $j=(\ell-1)N_{\ell}+i$, $i=1,\ldots,N_{\ell}$.

\begin{figure}[!b] % Figure 7, contours TTT
    \centering
    \begin{minipage}{0.9\textwidth} % container for the 2x2 block
        \centering
        % First row
        \includegraphics[clip,trim=7.2cm 12.1cm 7.5cm 9.5cm,width=0.46\textwidth]{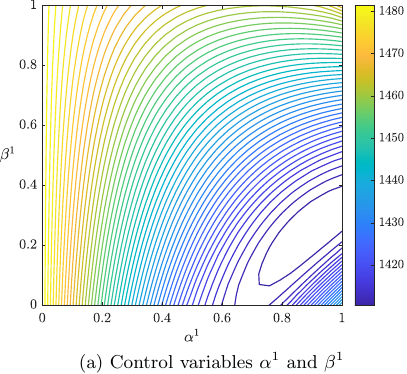}%
        \hspace{0.04\textwidth}%
        \includegraphics[clip,trim=7.2cm 12.1cm 7.5cm 9.5cm,width=0.46\textwidth]{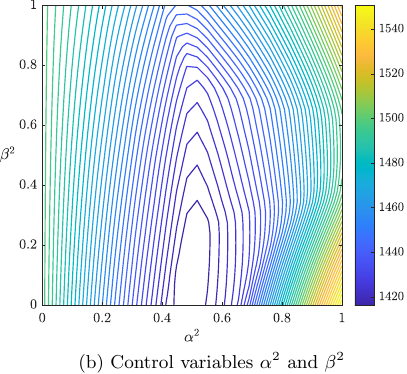}

        \vspace{0.7em} % space between rows

        % Second row
        \includegraphics[clip,trim=7.2cm 11.9cm 7.5cm 9.5cm,width=0.46\textwidth]{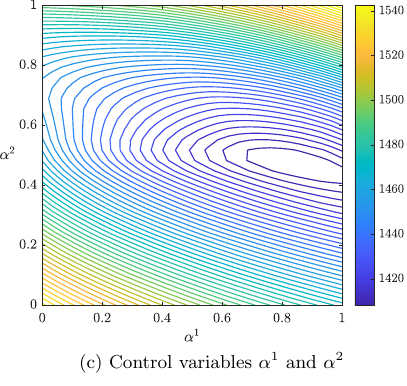}%
        \hspace{0.04\textwidth}%
        \includegraphics[clip,trim=7.2cm 11.9cm 7.5cm 9.5cm,width=0.46\textwidth]{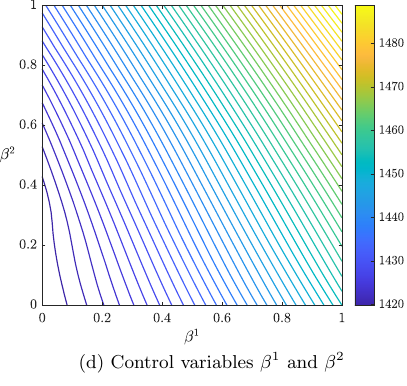}
    \end{minipage}

    \caption{Plot of the level curves for the TTT for two different control variables, where the fixed parameters have the values defined in \eqref{eq:num_tests:fixed_parameters}. In (a) and (b), the class is fixed, while in (c) and (d), the junction is fixed.}
    \label{fig:num_tests:func_val_TTT}
\end{figure}

The state variables for each time iteration $\nu$ are defined as $y^{\nu}=(y^{1,\nu},\ldots,y^{N,\nu})$, where $y^{c,\nu}=(l^{c,\nu},\rho^{c,\nu})$, $c=1,2$, are vectors of size $MT$, with $M=1+\sum_{i=1}^7 N_{\ell}$. The structure of the matrix $\partial_y E\in\mathcal{M}_{NMT\times NMT}(\R)$ is as presented in Figure \ref{fig:DAG:struct_matrix_E_y}, where, in particular,
\begin{equation*}
    \dfrac{\partial E^{c,\nu}}{\partial y^{c,\nu}} = I_M, \qquad \dfrac{\partial E^{c,\nu}}{\partial y^{g,\nu}} = 0 ~\text{ for }~ g\neq c,
\end{equation*}
and
\begin{equation*}
    \dfrac{\partial E^{c,\nu}}{\partial y^{g,\nu-1}} =
    \begin{pmatrix}
        \frac{\partial E_0^{c,\nu}}{\partial l^{g,\nu-1}} & \frac{\partial E_0^{c,\nu}}{\partial \rho_1^{g,\nu-1}} & 0 & 0 & \cdots & 0\\[8pt]
        \frac{\partial E_1^{c,\nu}}{\partial l^{g,\nu-1}} & \frac{\partial E_1^{c,\nu}}{\partial \rho_1^{g,\nu-1}} & \frac{\partial E_1^{c,\nu}}{\partial \rho_2^{g,\nu-1}} & 0 & \cdots & 0\\[8pt]
        0 & \frac{\partial E_2^{c,\nu}}{\partial \rho_1^{g,\nu-1}} & \frac{\partial E_2^{c,\nu}}{\partial \rho_2^{g,\nu-1}} & \frac{\partial E_2^{c,\nu}}{\partial \rho_3^{g,\nu-1}} & \cdots & 0\\
        \vdots & \vdots & \vdots & \vdots & \ddots & \vdots\\
        0 & 0  & 0 & 0 & \cdots & \frac{\partial E_{N_1+\cdots+N_7}^{c,\nu}}{\partial \rho_{N_1+\cdots+N_7}^{g,\nu-1}}
    \end{pmatrix}.
\end{equation*}
As mentioned above, the flow is controlled at junctions $e=2$ and $e=3$, where the control variables are the divergence parameters $\alpha_i^c$, shown in Figure \ref{fig:num_tests:example:network}. In particular, we now set them as
\begin{equation}
    \alpha_2^c = \alpha^c, \quad \alpha_4^c = \beta^c \quad (\text{and } \alpha_3^c = 1-\alpha^c, \quad \alpha_5^c = 1-\beta^c),
    \label{eq:num_tests:example:2_classes_4_control_points}
\end{equation}
with $0\leq \alpha^c, \beta^c\leq 1$ for every class $c$, for $c=1,2$.

\begin{figure}[!b] % Figure 8, contours TTD
    \centering
    \begin{minipage}{0.9\textwidth} % container for the 2x2 block
        \centering
        % First row
        \includegraphics[clip,trim=7.2cm 12.1cm 7.5cm 9.5cm,width=0.46\textwidth]{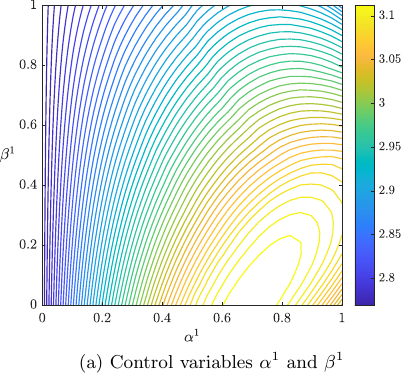}%
        \hspace{0.04\textwidth}%
        \includegraphics[clip,trim=7.2cm 12.1cm 7.5cm 9.5cm,width=0.46\textwidth]{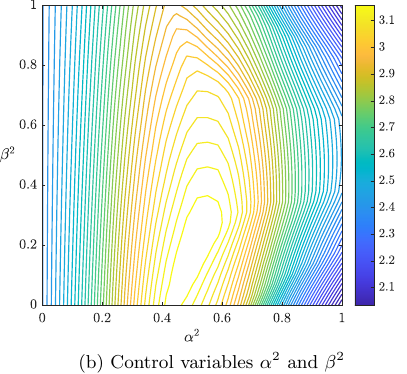}

        \vspace{0.7em} % space between rows

        % Second row
        \includegraphics[clip,trim=7.2cm 11.9cm 7.5cm 9.5cm,width=0.46\textwidth]{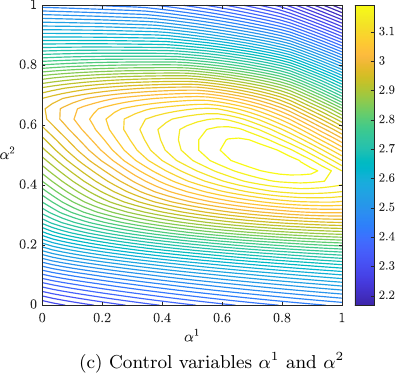}%
        \hspace{0.04\textwidth}%
        \includegraphics[clip,trim=7.2cm 11.9cm 7.5cm 9.5cm,width=0.46\textwidth]{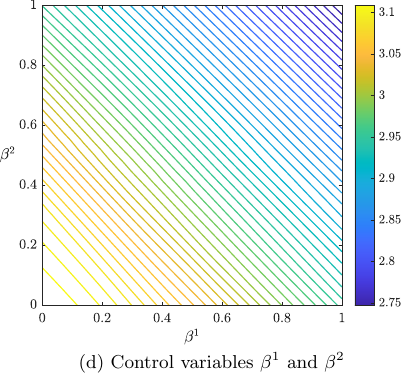}
    \end{minipage}

    \caption{Plot of the level curves for the TTD for two different control variables, where the fixed parameters have the values defined in \eqref{eq:num_tests:fixed_parameters}. In (a) and (b), the class is fixed, while in (c) and (d), the junction is fixed.}
    \label{fig:num_tests:func_val_TTD}
\end{figure}

For the first numerical test, taking the control points $\alpha^c$ and $\beta^c$, $c=1,2$, to be constant over the time interval $[0,T_f]$, we plot the two objective functions considered in Section \ref{sec_art:cost} for two out of the four control variables in \eqref{eq:num_tests:example:2_classes_4_control_points}, where the other two parameters are considered to be fixed, with values for $\alpha^c$ and $\beta^c$, $c=1,2$, when fixed, given as
\begin{equation}
    \alpha^1=1/2, \qquad \alpha^2=1/3, \qquad \beta^1=1/2, \qquad \beta^2=1/3. \label{eq:num_tests:fixed_parameters}
\end{equation}

\begin{table}[t] % Table 1
    \caption{Minimum values obtained for $J=\mathrm{TTT}$ and $J=\mathrm{TTD}$ depending on the class.}
    \label{table:num_tests:TTT_TTD_depending_on_class}
    \centering
    \begin{tabular}{lccc ccc}
        \toprule
                       & $\alpha^1$ & $\beta^1$   & $\min J$   & $\alpha^2$ & $\beta^2$ & $\min J$ \\
        \midrule
        $\mathrm{TTT}$ & 0.9310     & 0.2759      & 1409.0587  & 0.4828     & 0.1379    & 1413.6021 \\
        $\mathrm{TTD}$ & 0          & $b\in[0,1]$ & 27623.8203 & 1          & 0         & 20107.7678 \\
        \bottomrule
    \end{tabular}
\end{table}

\begin{table}[t] % Table 2
    \caption{Minimum values obtained for $J=\mathrm{TTT}$ and $J=\mathrm{TTD}$ depending on the junction.}
    \label{table:num_tests:TTT_TTD_depending_on_junction}
    \centering
    \begin{tabular}{lccc ccc}
        \toprule
                       & $\alpha^1$ & $\alpha^2$ & $\min J$   & $\beta^1$ & $\beta^2$ & $\min J$ \\
        \midrule
        $\mathrm{TTT}$ & 0.8966     & 0.4483     & 1405.4970  & 0         & 0         & 1418.4881 \\
        $\mathrm{TTD}$ & 1          & 1          & 21472.1444 & 1         & 1         & 27402.6890 \\
        \bottomrule
    \end{tabular}
\end{table}

Figures \ref{fig:num_tests:func_val_TTT} and \ref{fig:num_tests:func_val_TTD} show the plots of the level curves for the total travel time \eqref{eq:TTT} and total travel distance \eqref{eq:TTD} functionals, depending on the class (\ref{fig:num_tests:func_val_TTT}a, \ref{fig:num_tests:func_val_TTT}b, \ref{fig:num_tests:func_val_TTD}a, \ref{fig:num_tests:func_val_TTD}b) and on the junction (\ref{fig:num_tests:func_val_TTT}c, \ref{fig:num_tests:func_val_TTT}d, \ref{fig:num_tests:func_val_TTD}c, \ref{fig:num_tests:func_val_TTD}d), respectively.

Tables \ref{table:num_tests:TTT_TTD_depending_on_class} and \ref{table:num_tests:TTT_TTD_depending_on_junction} present the minimum values obtained for the TTT and TTD in Figures \ref{fig:num_tests:func_val_TTT} and \ref{fig:num_tests:func_val_TTD}, as well as the respective values of the different control parameters where such minimum is reached. These results show that the model identifies a different optimal performance under two control points for both functionals. In the following, we will consider the four control points in \eqref{eq:num_tests:example:2_classes_4_control_points} as time-dependent to evaluate the benefits of a dynamic control strategy.

\subsection{Gradient computation: discrete adjoint VS finite differences}
\label{sec_art:num_tests:comparison_DAG_finitediff}

This test aims to show that, for a sufficiently large number of control points, adjoint gradient computation becomes competitive with respect to finite differences. To this end, we consider time-dependent controls: given $n\geq1$, we divide the simulation horizon in subintervals as
\begin{equation*}
    [0,T_f[ = \bigcup_{j=1}^{n}\left[\frac{j-1}{n}T_f,\frac{j}{n}T_f\right[,
\end{equation*}
and define the control actions $u_j=(u_j^1,u_j^2)$, with $u^c_j=(\alpha^c_j,\beta^c_j)$, $c=1,2$, to be piece-wise constant on each subinterval $\left[\frac{j-1}{n}T_f,\frac{j}{n}T_f\right[$, $j=1,\ldots,n$. Following \eqref{eq:num_tests:example:2_classes_4_control_points}, the matrix $\partial_u E$ can be reduced to $\partial_u E\in\mathcal{M}_{NMT\times 4n}(\mathbb{R})$, with $u=(u_1,\ldots,u_n)$, where $\partial_u E = (\partial_{u_j} E^{\nu})_{\nu=1,\ldots,T,\;j=1,\ldots,n}$ and
\begin{equation*}
    \dfrac{\partial E^{\nu}}{\partial u_j} =
    \left(
    \renewcommand{\arraystretch}{1.5}
    \begin{array}{c|c}
        \frac{\partial E^{1,\nu}}{\partial u^1_j} & 0\\[3pt] \hline
        0 & \frac{\partial E^{2,\nu}}{\partial u^2_j}
    \end{array}
    \right)\in\mathcal{M}_{NM\times 4}(\mathbb{R}), \quad \dfrac{\partial E^{c,\nu}}{\partial u^c_j} = 
    \left( 
    \renewcommand{\arraystretch}{1.5}
    \begin{array}{c|c}
        \dfrac{\partial E^{c,\nu}}{\partial \alpha^c_j} & \dfrac{\partial E^{c,\nu}}{\partial \beta^c_j}
    \end{array}
    \right)\in\mathcal{M}_{M\times 2}(\mathbb{R}).
\end{equation*}
To solve the minimization problem \eqref{eq:OptProbl} for the total travel time~\eqref{eq:TTT}, we used the nonlinear optimization routine \texttt{fmincon} available in the \textsc{Matlab} R2024a library, with \texttt{sqp} as the chosen optimization algorithm. The maximum number of function evaluations (fval) and iterations (iter) are set to $10^4$ and $10^3$, respectively. The optimality and constrain tolerances, denoted as $\text{TOL}_{\text{opt}}$ and $\text{TOL}_{\text{cons}}$, are given a value of $10^{-6}$, while the function and step tolerances, i.e., $\text{TOL}_{\text{fun}}$ and $\text{TOL}_{\text{step}}$, are $10^{-7}$. The initial values given to the control variables $\alpha^c$ and $\beta^c$, $c=1,2$, are 
\begin{equation}
    \alpha^1 = 0.7, \qquad \alpha^2 = 0.4, \qquad \beta^1 = 0.1, \qquad \beta^2 = 0.7. \label{eq:num_tests:comparison:initial_cond}
\end{equation}
The goal is to compare the outcome of \eqref{eq:OptProbl} when the gradient is calculated through central finite differences with respect to  the discrete adjoint gradient method described in Section~\ref{sec_art:DAG}. 

As the definition of the total travel time \eqref{eq:TTT} does not depend explicitly on the control parameters $u$, then $\partial_u \text{TTT}=0$, which reduces the cost gradient \eqref{eq:costgradient} to
\begin{equation*}
    \dfrac{d}{du} \text{TTT}(u,y(u)) = \xi^T\,\partial_u E,
\end{equation*}
where $\xi$ solves the adjoint equation \eqref{eq:adjointeq}.

We set the number of time subintervals as $n=2^a$, $a\in\mathbb{N}_{\geq 0}$, i.e. there are $2^{a+2}$ control variables. For the numerical tests, the chosen values of $a=\log_2 n$ are $a=0,\ldots,6$.

\begin{figure}[!ht] % Figure 9
    \centering
    % Entire 2x2 block centered as a single entity
    \begin{minipage}{0.96\textwidth} % adjust total width to fit images + spacing
        \centering
        % First row
        \includegraphics[width=0.48\textwidth]{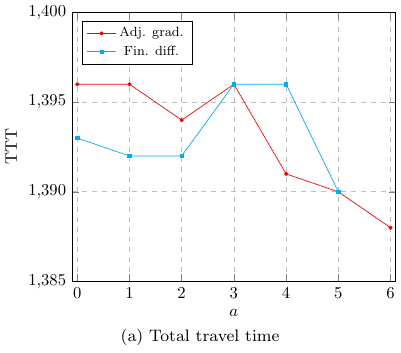}%
        \hspace{0.04\textwidth}%
        \includegraphics[width=0.46\textwidth]{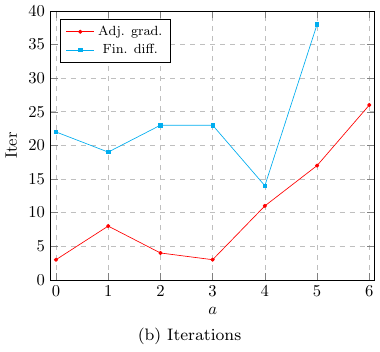}

        \vspace{0.5em} % space between rows

        % Second row
        \includegraphics[width=0.46\textwidth]{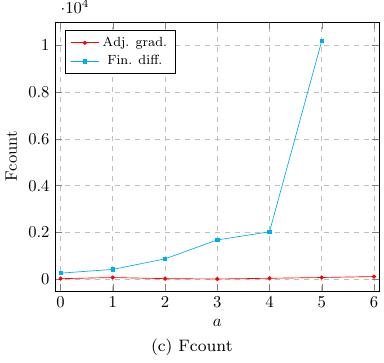}%
        \hspace{0.03\textwidth}%
        \includegraphics[width=0.48\textwidth]{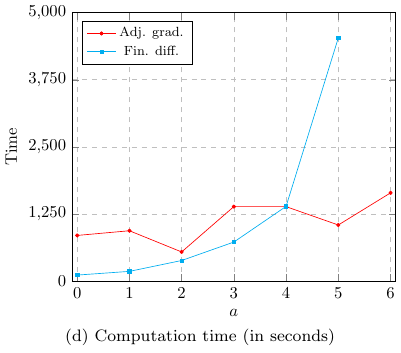}
    \end{minipage}

    \caption{Results comparing the discrete adjoint gradient and the finite difference computation for a \textsc{Matlab}~R2024a \texttt{fmincon} optimization with an \texttt{sqp} algorithm, with increasing values of $a=\log_2 n$, $a=0,\ldots,6$, where $n$ is the number of time-interval subdivisions.}
    \label{fig:num_tests:results_plot_a}
\end{figure}

Figure \ref{fig:num_tests:results_plot_a} compares the performances of the discrete adjoint gradient method and the finite difference method in solving the optimization problem as the number of time-interval subdivisions $n$ increases. In particular:
\begin{itemize}
    \item Figure \ref{fig:num_tests:results_plot_a}a displays the total travel time obtained for different values of $a$, where, for a higher number of time-interval subdivisions, the adjoint gradient provides in general a smaller value for the TTT compared to the finite differences.

    \item In Figure \ref{fig:num_tests:results_plot_a}b we can see that the number of iterations required for convergence is smaller for the adjoint gradient that for the finite differences for every value of $a$;

    \item Figure \ref{fig:num_tests:results_plot_a}c presents the function evaluation count, where the adjoint gradient requires a significantly smaller number of function evaluations compared to finite differences;
    
    \item the computational time in seconds required by each method is shown in Figure \ref{fig:num_tests:results_plot_a}d, where, for larger values of $n$, the adjoint gradient method is computationally more efficient than finite differences.
\end{itemize}

The results show that the adjoint gradient method outperformas finite differences as the number of control points grows, while for small $n$ the finite difference method yields better results. Although the finite-difference method may provide more optimal results for smaller values of $n$, it experiences a rapid increase in the number of iterations, function evaluations and overall computation time, which confirms it as undesirable method when considering a higher number of time-interval subdivisions. Moreover, for $a>5$, the optimization process stops before obtaining an optimal result for the finite difference method. On the other hand, the discrete adjoint gradient maintains a stable number of function evaluations and a small increase in computation time for any value of $n$.  

Therefore, in cases where a high number of time-dependent control points is required, the adjoint gradient is the more reliable and computationally efficient choice, as the exponential growth for the optimization considering the finite difference method is inefficient for large-scale problems.

These results have practical implications for real-time and large-scale traffic management. Modern control systems rely on continuously updated data provided by sensors and GPS-equipped vehicles. Navigation systems, for example, dynamically adjust routes based on congestion levels and travel-time estimates, requiring a repeated and rapid optimization of control parameters.

Consequently, the proposed adjoint-based optimization can be naturally integrated into real-time traffic control systems. In particular, the method can efficiently compute class-specific control adjustments, such as rerouting heavy vehicles or prioritizing faster lanes, within rapidly changing traffic environments, contributing directly to relieving congestion and improving network performance.

\subsection{Pareto fronts}
\label{sec_art:num_tests:Pareto_front}

\begin{figure}[!ht] % Figure 10
	\centering
	\includegraphics[scale=1.05]{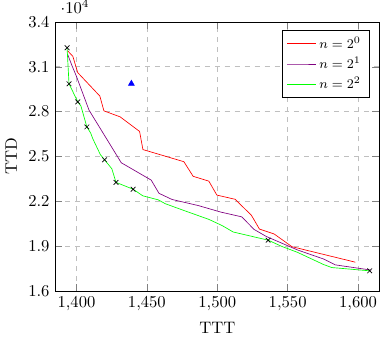} % Pareto Front (TTT,TTD)
    \caption[]{Plot of the Pareto front of the TTT against the TTD for different time-interval subdivisions $n$.}
    \label{fig:num_tests:ParetoFront}
\end{figure}

For better traffic management, it is also interesting to consider a multi-objective optimization approach to illustrate the potential of the proposed framework. In this setting, the total travel time \eqref{eq:TTT} and the total travel distance \eqref{eq:TTD} represent conflicting objectives: minimizing the TTT promotes faster circulation and reduced delays, while minimizing the TTD favors shorter routes, which may lead to increased congestion and longer travel times. To capture this trade-off, we compute the Pareto front, which represents the set of non-dominated optimal solutions, each corresponding to a different balance between these two functionals. The Pareto fronts are here computed using the \textsc{Matlab} routine \texttt{gamultiobj}, applied to the same network and control settings as in the previous subsection. 

Figure \ref{fig:num_tests:ParetoFront} presents a Pareto front plot where each plotted curve corresponds to a different number of time-interval subdivisions $n$, and the node represented by a triangle denotes the TTT and TTD obtained for the initial $\alpha^c$ and $\beta^c$ values defined in \eqref{eq:num_tests:comparison:initial_cond}, i.e. $\text{TTT}=1439$ and $\text{TTD}=29865$. Moreover, if we want to reduce both functionals compared to such initial point, the chosen solution must lie within the region bounded by the projections of the point onto the $x$ and $y$-axis. Points lying “between'' these projections correspond to configurations for which both TTT and TTD decrease, reflecting an overall improvement relative to the starting control values.

When the value of $n$ increases, the higher number of control points leads to a decrease in the ranks for the solutions of such $n$, enabling more refined adjustments of traffic over time and more control flexibility. This typically results in solutions that achieve a better balance between TTT and TTD.

\begin{figure}[!ht] % Figure 11
    \centering
    \begin{minipage}{0.96\textwidth} % container for the 2x2 block
        \centering
        % First row
        \includegraphics[clip,trim=7.3cm 17.5cm 7.5cm 4.7cm,width=0.49\textwidth]{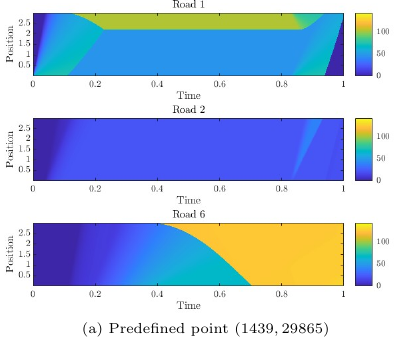}%
        \hspace{0.019\textwidth}%
        \includegraphics[clip,trim=7.3cm 17.5cm 7.5cm 4.7cm,width=0.49\textwidth]{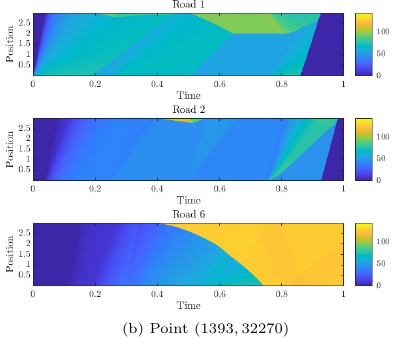}

        \vspace{0.8em} % space between rows

        % Second row
        \includegraphics[clip,trim=7.3cm 17.3cm 7.5cm 4.7cm,width=0.49\textwidth]{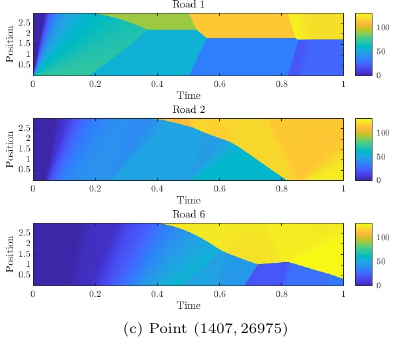}%
        \hspace{0.019\textwidth}%
        \includegraphics[clip,trim=7.3cm 17.3cm 7.5cm 4.7cm,width=0.49\textwidth]{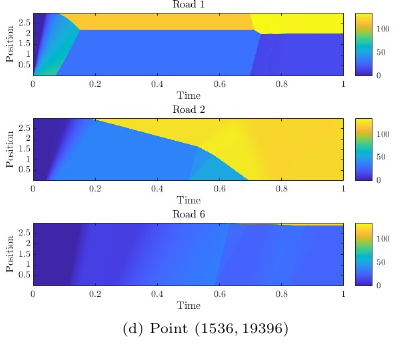}
    \end{minipage}

    \caption{Results comparing the total vehicle densities on roads 1, 2 and 6 of the network represented in Figure \ref{fig:num_tests:example:network} for some of the points shown in Figure \ref{fig:num_tests:ParetoFront}, together with the initial point $(1439,29865)$. Point (b) minimizes TTT over TTD, point (c) minimizes both functionals and point (d) minimizes TTD over TTT.}
    \label{fig:num_tests:PF:several_points}
\end{figure}

From a practical point of view, the Pareto front provides valuable insight for decision-making in traffic management. Points located towards the lower-left corner of the plot (smaller TTT and TTD) correspond to balanced strategies that reduce both congestion and travelled distance. In contrast, points with low TTT but higher TTD correspond to control actions prioritizing faster travel even if the chosen route is longer, while points with low TTD and higher TTT represent strategies minimizing total mileage and route length, but at the cost of slower movement and potential congestion. This trade-off analysis can therefore assist transportation planners or automated systems in selecting control strategies that align with specific operational priorities.

Figure \ref{fig:num_tests:PF:several_points} shows the time evolution of the total vehicular density on selected roads for different Pareto-optimal solutions. Only roads 1, 2, and 6 of the considered network are displayed, since these links exhibit the most significant variations in density across all the chosen Pareto points. In contrast, roads 3, 4, 5, and 7 show only marginal changes with respect to the initial configuration, and their inclusion would not provide additional insight into the effects of the optimization. The density changes observed on roads 1, 2, and 6 indicate that the optimization primarily redistributes traffic on these links. More precisely, the effects of the optimization can be observed differently on each road.

On road 1, Pareto-optimal solutions prioritizing the minimization of TTT over TTD are associated with a reduction of the total density and a smoother temporal evolution, indicating a more efficient exit of vehicles from this link. On road 2, density profiles corresponding to TTD-oriented solutions exhibit higher density levels, reflecting the tendency to concentrate traffic along shorter routes. Finally, road 6 appears to be particularly sensitive to the chosen Pareto point: balanced solutions moderate density peaks while maintaining temporal regularity, whereas more extreme trade-off choices produce sharper density accumulations or a faster exit, depending on the prioritized objective.

\begin{table}[!ht] % Table 3
    \centering
    \caption{Percentage change in TTT and TTD for several Pareto-optimal solutions for $n=4$ in Figure \ref{fig:num_tests:ParetoFront} compared to the initial point of the optimization. The $+$ and $-$ signs denote a reduction or an increase of the value of the considered functional, respectively.}
    \label{table:num_tests:Pareto_Front_improvement}
    \begin{tabular}{lll} 
         \toprule
         \text{Points} & $\%\Delta\text{TTT}$ & $\%\Delta\text{TTD}$ \\
         \midrule
         (1393,32270) & $+3.20\%$  & $-8.05\%$ \\
         (1395,29860) & $+3.06\%$  & $+0.02\%$ \\
         (1401,28650) & $+2.64\%$  & $+4.07\%$ \\
         (1407,26975) & $+2.22\%$  & $+9.68\%$ \\
         (1420,24776) & $+1.34\%$  & $+17.04\%$ \\
         (1428,23262) & $+0.76\%$  & $+22.11\%$ \\
         (1440,22805) & $-0.07\%$  & $+23.64\%$ \\
         (1536,19396) & $-6.74\%$  & $+35.07\%$ \\
         (1608,17344) & $-11.74\%$ & $+41.93\%$ \\
         \bottomrule
    \end{tabular} 
\end{table}

Table \ref{table:num_tests:Pareto_Front_improvement} illustrates the practical trade-off captured by the Pareto front by comparing several optimal solutions for $n=4$ against the initial optimization point, marked in Figure \ref{fig:num_tests:ParetoFront} by $\times$. The columns $\%\Delta\text{TTT}$ and $\%\Delta\text{TTD}$ represent the percentage change in total travel time and total travel distance, respectively, compared to the initial value given by \eqref{eq:num_tests:comparison:initial_cond}. A positive value indicates an improvement (reduction) with respect to the original value, while a negative percentage represents a deterioration (increase). This interpretation is consistent with the geometric observation above: points that lie between the projections of the initial triangular marker on the Pareto plane are precisely those for which both $\%\Delta\text{TTT}$ and $\%\Delta\text{TTD}$ are positive, indicating a reduction in travel time and travelled distance. Points falling outside this region show the expected trade-off behaviour, with one metric improving at the expense of the other. For instance, the point $(1407,26975)$ represents a balanced strategy, achieving a $2.22\%$ reduction in TTT and a $9.68\%$ in TTD. In contrast, $(1608,17344)$ shows that the minimization of distance is prioritized over the minimization of time, as there is a significant $41.93\%$ reduction of TTD but an increase of $11.74\%$ for the TTT.

However, we observed that increasing the number of time-interval subdivisions considered beyond $n=4$ deteriorates the precision of the algorithm. A moderate increase of $n$ allows for more refined results, but exceeding an optimal threshold reduces the accuracy and consistency of the results, making the algorithm unable to converge to the expected solution.

This decrease in accuracy highlights the need for more performant Pareto front computation algorithms and an optimal selection of $n$. While more control points provide flexibility for adjusting traffic flow, taking an excessive number of subdivisions of $[0,T_f]$ can lead to over-parametrization, resulting in ill-conditioned optimization problems and an increase in numerical errors.

\section{Conclusions}
\label{sec_art:conclusions}

This article presented a discrete adjoint gradient computation method applied to a multi-class macroscopic traffic flow model on road networks. The model was formulated as a system of hyperbolic conservation laws, and was discretized using a Godunov-type finite volume scheme based on demand and supply functions, taking into account the traffic dynamics at junctions. The adjoint-based optimization framework was then employed to efficiently compute gradients to be used in optimization routines.

Numerical tests demonstrated the effectiveness of the proposed approach in optimizing the total travel time and total travel distance, important traffic performance metrics. The comparison between adjoint gradient computation and finite differences showed that, while finite differences may perform well for a small number of control points, the adjoint method remains computationally efficient as the number of control parameters increases. Moreover, the computational efficiency of the discrete adjoint gradient approach makes it suitable for integration into real-time traffic management settings and GPS-based routing systems, where optimization must be repeated within short time periods for the changing traffic conditions.

Additionally, the multi-objective optimization using Pareto front analysis showed the relationship between these different traffic performance metrics and the need for an optimal selection of number of time-dependent control points when considering several metrics.

Future work could explore the extension of this method to larger and more complex networks, incorporating additional control variables and other cost functionals such as fuel consumption, network throughput or $\text{CO}_2$ emissions, as well as integrating real-time traffic data or adaptive control strategies.

\backmatter

\section*{Statements and Declarations}

\subsubsection*{Acknowledgements}
The authors would like to thank the members of the Datahyking Doctoral Network for their valuable discussions and support throughout this research.

\subsubsection*{Funding} 
This work was funded by the European Union’s Horizon Europe research and innovation programme under the Marie Sklodowska-Curie Doctoral Network Datahyking (Grant No. 101072546).

\subsubsection*{Competing Interests} 
The authors declare that they have no competing interests.

\subsubsection*{Ethics Approval and Consent to Participate} 
This study does not involve human participants or animals and therefore does not require ethics approval or consent to participate.

\subsubsection*{Consent for Publication} 
Not applicable.

\subsubsection*{Data Availability} 
The data that support the findings of this study are available from the corresponding author, C.~Mezquita-Nieto, upon reasonable request. 

\subsubsection*{Code Availability} 
The numerical methods and algorithms used in this study are described in detail in the manuscript. The code implementation used to produce the numerical results is available from the corresponding author upon reasonable request.

\subsubsection*{Author Contributions}
\begin{itemize}
    \item[-] Conceptualization: P. Goatin, C. Mezquita-Nieto.
    \item[-] Methodology: P. Goatin, A. Klar, C. Mezquita-Nieto.
    \item[-] Software: C. Mezquita-Nieto.
    \item[-] Formal analysis: P. Goatin, C. Mezquita-Nieto.
    \item[-] Validation: C. Mezquita-Nieto.
    \item[-] Writing -- original draft: C. Mezquita-Nieto.
    \item[-] Writing -- review and editing: P. Goatin, A. Klar, C. Mezquita-Nieto.
\end{itemize}
All authors read and approved the final manuscript.

\begin{appendices}

\section{Details on DAG computation}
\label{app:DAG}

\subsection{Update equation}
\label{app:DAG:E}

Denote $s$ as a variable that can be either a state or a control variable, $s\in u\cup y$. Considering a particular update equation $E_{\ell,j}^{c,\nu}$ \eqref{eq:DiscModelEq_rho}-\eqref{eq:DiscModelEq_l} and applying the chain rule, its partial derivatives with respect to the state variables $y$ are
\begin{align*}
    \partial_s E_{\ell,j}^{c,\nu} &= 
    \begin{cases}
        1 & \text{if } s = \rho_{\ell,j}^{c,\nu},\, j=1,\ldots,N_{\ell},\\[8pt]
        -1 + \lambda \left[ \partial_s F_{\ell,j+1/2}^{c,\nu-1} - \partial_s F_{\ell,j-1/2}^{c,\nu-1} \right] & \text{if } s = \rho_{\ell,j}^{c,\nu-1}, \, j=2,\ldots,N_{\ell}-1,\\[8pt]
        \lambda \left[ \partial_s F_{\ell,j+1/2}^{c,\nu-1} - \partial_s F_{\ell,j-1/2}^{c,\nu-1} \right] & \text{if } s = \rho_{\ell,j}^{g,\nu-1}, \, j=2,\ldots,N_{\ell}-1,\, g\neq c,\\[8pt]
        \lambda\,\partial_s F_{\ell,j+1/2}^{c,\nu-1} & \text{if } s = \rho_{\ell,j+1}^{g,\nu-1},\, j=2,\ldots,N_{\ell}-1,\, \forall g,\\[8pt]
        - \lambda\,\partial_s F_{\ell,j-1/2}^{c,\nu-1} & \text{if } s = \rho_{\ell,j-1}^{g,\nu-1},\, j=2,\ldots,N_{\ell}-1,\, \forall g,\\[8pt]
        0 & \text{else,}
    \end{cases}\\[5pt]
    \partial_s E_{\ell,0}^{c,\nu} &= 
    \begin{cases}
        1 & \text{if } s = l_{\ell}^{c,\nu}, \\[5pt]
        -1 & \text{if } s = l_{\ell}^{c,\nu-1}, \, l_{\ell}^{c,\nu-1} + \Delta t (F_{\ell,in}^c-\hat{\gamma}^{c,\nu-1}_{\ell,in}) > 0,\\[5pt]
        \Delta t\, \partial_s \hat{\gamma}_{\ell,in}^{c,\nu-1} & \text{if } s=\rho_{\ell,1}^{g,\nu-1}, \, l_{\ell}^{c,\nu-1} + \Delta t (F_{\ell,in}^c-\hat{\gamma}^{c,\nu-1}_{\ell,in}) > 0,\, \forall g,\\[5pt]
        0 & \text{else,}
        \end{cases}
\end{align*}
while for the control variables $u$ they are
\begin{align*}
    \partial_s E_{\ell,j}^{c,\nu} &= 
    \begin{cases}
        \lambda \left[ \partial_s F_{\ell,j+1/2}^{c,\nu-1} - \partial_s F_{\ell,j-1/2}^{c,\nu-1} \right] & \text{if } s\in u,\, j=1,\ldots,N_{\ell},\\
        0 & \text{else,}
    \end{cases}\\[5pt]
    \partial_s E_{\ell,0}^{c,\nu} &= 
    \begin{cases}
        \Delta t \, \partial_s \hat{\gamma}_{\ell,in}^{c,\nu-1} & \text{if } s\in u, \, l_{\ell}^{c,\nu-1} + \Delta t (F_{\ell,in}^c-\hat{\gamma}^{c,\nu-1}_{\ell,in}) > 0,\\
        0 & \text{else.}
    \end{cases}
\end{align*}

\subsection{Godunov scheme}
\label{app:DAG:Godunov}

Returning to the finite volume Godunov scheme \eqref{eq:FV}, its partial derivatives with respect to the state variables are
\begin{equation*}
    \partial_s F_{\ell,j+1/2}^{c,\nu} =
    \begin{cases}
        k_{\ell,j}^{c,\nu} \, \partial_s S_{\ell}^c(r_{\ell,j+1}^{\nu}) & \text{if } s=\rho_{\ell,j+1}^{c,\nu},\, \ell\in\mathcal{S}^{\text{G}}_2,\, \forall g,\\[8pt]
        (k_{\ell,j}^{c,\nu})^{'} D_{\ell}^c(r_{\ell,j}^{\nu}) + k_{\ell,j}^{c,\nu} \,\partial_s D_{\ell}^c(r_{\ell,j}^{\nu}) & \text{if } s=\rho_{\ell,j}^{c,\nu},\, \ell\in\mathcal{S}^{\text{G}}_1,\\[9pt]
        k_{\ell,j}^{c,\nu}\left[-\dfrac{D^c_{\ell}(r_{\ell,j}^{\nu})}{r_{\ell,j}^{\nu}} + \partial_s D^c_{\ell}(r_{\ell,j}^{\nu})\right] & \text{if } s=\rho_{\ell,j}^{g,\nu},\, \ell\in\mathcal{S}^{\text{G}}_1,\, g\neq c,\\[13pt]
        (k_{\ell,j}^{c,\nu})^{'} S_{\ell}^c(r_{\ell,j+1}^{\nu}) & \text{if } s=\rho_{\ell,j}^{c,\nu},\, \ell\in\mathcal{S}^{\text{G}}_2,\\[8pt]
        -\dfrac{k_{\ell,j}^{c,\nu}}{r_{\ell,j}^{\nu}} \, S_{\ell}^c(r_{\ell,j+1}^{\nu}) & \text{if } s=\rho_{\ell,j}^{g,\nu},\, \ell\in\mathcal{S}^{\text{G}}_2,\, g\neq c,\\
        0 & \text{else,}
    \end{cases}
\end{equation*}
where $k_{\ell,j}^{c,\nu}$ and $(k_{\ell,j}^{c,\nu})'$ are defined as
\begin{equation}
    k_{\ell,j}^{c,\nu}\coloneqq \dfrac{\rho_{\ell,j}^{c,\nu}}{r_{\ell,j}^{\nu}} \quad \text{and} \quad (k_{\ell,j}^{c,\nu})^{'}\coloneqq \dfrac{\partial k_{\ell,j}^{c,\nu}}{\partial \rho_{\ell,j}^{c,\nu}} = \dfrac{r_{\ell,j}^{\nu}-\rho_{\ell,j}^{c,\nu}}{(r_{\ell,j}^{\nu})^2}, \label{eq:k_quotients}
\end{equation}
and
\begin{align*}
    \mathcal{S}_1^{\text{G}} &\coloneqq \big\{ \ell\; :\; D_{\ell}^c(r_{{\ell},j}^{\nu})\leq S_{\ell}^c(r_{{\ell},j+1}^{\nu}) \big\},\\
    \mathcal{S}_2^{\text{G}} &\coloneqq \big\{ \ell\; :\; D_{\ell}^c(r_{{\ell},j}^{\nu})> S_{\ell}^c(r_{{\ell},j+1}^{\nu}) \big\},\\
    \mathcal{S}_1^{\text{G}}&\cap\mathcal{S}_2^{\text{G}} = \varnothing.
\end{align*}
Similarly, the derivatives with respect to the control variables are
\begin{equation*}
    \partial_s F_{\ell,j+1/2}^{c,\nu} =
    \begin{cases}
        k_{\ell,j}^{c,\nu} \, \partial_s D_{\ell}^c(r_{\ell,j}^{\nu}) & \text{if } s=V_{\ell,c},\, \ell\in\mathcal{S}^{\text{G}}_1, \\[3pt]
        k_{\ell,j}^{c,\nu} \, \partial_s S_{\ell}^c(r_{\ell,j+1}^{\nu}) & \text{if } s=V_{\ell,c},\, \ell\in\mathcal{S}^{\text{G}}_2,\\
        0 & \text{else.}
    \end{cases}
\end{equation*}

\subsection{Demand and supply}
\label{app:DAG:DS}

The partial derivatives with respect to $s$ of the demand and the supply functions are
\begin{align*}
    \partial_s D_{\ell}^c(r_{\ell,j}^{\nu}) &=
    \begin{cases}
        v_{\ell,c}(r_j^{\nu}) + r_{\ell,j}^{\nu}\, \partial_s v_{\ell,c}(r_{\ell,j}^{\nu}) & \text{if } s=\rho_{\ell,j}^{g,\nu}, \, r_{\ell,j}^{\nu}\leq r_{\ell,j,cr}^{c,\nu},\, \forall g,\\[3pt]
        \min \lbrace r_{\ell,j}^{\nu},r_{\ell,j,cr}^{c,\nu} \rbrace\, \partial_sv_{\ell,c}(\min \lbrace r_{\ell,j}^{\nu},r_{\ell,j,cr}^{c,\nu} \rbrace) & \text{if } s=V_{\ell,c},\\
        0 & \text{else,}
    \end{cases}\\[5pt]
    \partial_s S_{\ell}^c(r_{\ell,j}^{\nu}) &=
    \begin{cases}
        v_{\ell,c}(r_{\ell,j}^{\nu}) + r_{\ell,j}^{\nu} \, \partial_s v_{\ell,c}(r_{\ell,j}^{\nu}) & \text{if } s=\rho_{\ell,j}^{g,\nu},\, r_{\ell,j}^{\nu}> r_{\ell,j,cr}^{c,\nu},\, \forall g,\\[3pt]
        \max \lbrace r_{\ell,j}^{\nu},r_{\ell,j,cr}^{c,\nu} \rbrace \, \partial_s v_{\ell,c}(\max \lbrace r_{\ell,j}^{\nu},r_{\ell,j,cr}^{c,\nu} \rbrace) & \text{if } s=V_{\ell,c},\\
        0 & \text{else.}
    \ \end{cases}
\end{align*}

\subsection{\texorpdfstring{$M\times 1$}{Mx1} merging junction} \label{app:DAG:Mx1}

For junctions of the type $M\times 1$, the flux $\hat{\gamma}_{M+1,in}^{c,\nu}$ is defined as
\begin{equation*}
    \hat{\gamma}_{M+1,in}^{c,\nu} = \sum_{i=1}^M \hat{\gamma}_{i,out}^{c,\nu}, \qquad \hat{\gamma}_{i,out}^{c,\nu} = k_{i,N_i}^{c,\nu} \, \beta_1^i,
\end{equation*}
where
\begin{align*}
    \beta_1^i &\coloneqq \min \left\{ D^c_i(r_{i,N_i}^{\nu}),\, \beta_2^i \right\}, \\
    \beta_2^i &\coloneqq \max \big\{ p_i^c\, S^c_{M+1}(r_{M+1,1}^{\nu}), S^c_{M+1}(r_{M+1,1}^{\nu}) - \sum_{j\neq i} D^c_j(r_{j,N_j}^{\nu}) \big\},
\end{align*}
and $k_{i,N_i}^{c,\nu}$ follows \eqref{eq:k_quotients}. Let
\begin{align*}
    \mathcal{S}_1^{\text{m}} &\coloneqq \big\{ i\; :\; \beta_1^i = D_i^c(r_{i,N_i}^{\nu}) \big\},\\[3pt]
    \mathcal{S}_2^{\text{m}} &\coloneqq \big\{ i\; :\; \beta_1^i = \beta_2^i = p_i^c\, S_{M+1}^c(r_{M+1,1}^{\nu}) \big\},\\[3pt]
    \mathcal{S}_3^{\text{m}} &\coloneqq \big\{ i\; :\; \beta_1^i = \beta_2^i = S^c_{M+1}(r_{M+1,1}^{\nu}) - \sum_{j\neq i} D^c_j(r_{j,N_j}^{\nu}) \big\},\\
    \mathcal{S}_i^{\text{m}}&\cap\mathcal{S}_j^{\text{m}} = \varnothing \quad\forall j\neq i.
\end{align*}
The partial derivatives of $\hat{\gamma}_{M+1,in}^{c,\nu}$ with respect to $s$ are obtained as
\begin{equation*}
    \partial_s \hat{\gamma}_{M+1,in}^{c,\nu} = \sum_{i=1}^M \partial_s \hat{\gamma}_{i,out}^{c,\nu},
\end{equation*}
and the derivatives of the fluxes $\hat{\gamma}_{i,out}^{c,\nu}$ with respect to the state variables are, for $i=1,\ldots,M$,

\begin{equation*}
    \partial_s \hat{\gamma}_{i,out}^{c,\nu} =
    \begin{cases}
        (k^{c,\nu}_{i,N_i})^{'} D_i^c(r_{i,N_i}^{\nu}) + k^{c,\nu}_{i,N_i} \partial_s D_i^c(r_{i,N_i}^{\nu}) & \text{if } s=\rho_{i,N_i}^{c,\nu},\, i\in\mathcal{S}_1^{\text{m}},\\[9pt]
        k^{c,\nu}_{i,N_i}\left[-\dfrac{D_i^c(r_{i,N_i}^{\nu})}{r_{i,N_i}^{\nu}} + \partial_s D_i^c(r_{i,N_i}^{\nu})\right] & \text{if } s=\rho_{i,N_i}^{g,\nu},\, i\in\mathcal{S}_1^{\text{m}},\, g\neq c,\\[15pt]
        (k^{c,\nu}_{i,N_i})^{'} \, p_i^c\, S_{M+1}^c(r_{M+1,1}^{\nu}) & \text{if } s=\rho_{i,N_i}^{c,\nu},\, i\in\mathcal{S}_2^{\text{m}},\\[13pt]
        -\dfrac{k^{c,\nu}_{i,N_i}}{r_{i,N_i}^{\nu}}\, p_i^c\, S_{M+1}^c(r_{M+1,1}^{\nu}) & \text{if } s=\rho_{i,N_i}^{g,\nu},\, i\in\mathcal{S}_2^{\text{m}},\, g\neq c,\\[13pt]
        (k^{c,\nu}_{i,N_i})^{'} \big( S^c_{M+1}(r_{M+1,1}^{\nu}) - \displaystyle\sum_{j\neq i} D^c_j(r_{j,N_j}^{\nu}) \big) & \text{if } s=\rho_{i,N_i}^{c,\nu},\, i\in\mathcal{S}_3^{\text{m}},\\[13pt]
        -\dfrac{k^{c,\nu}_{i,N_i}}{r_{i,N_i}^{\nu}} \big( S^c_{M+1}(r_{M+1,1}^{\nu}) - \displaystyle\sum_{j\neq i} D^c_j(r_{j,N_j}^{\nu}) \big) & \text{if } s=\rho_{i,N_i}^{g,\nu},\, i\in\mathcal{S}_3^{\text{m}},\, g\neq c,\\[13pt]
        -k^{c,\nu}_{i,N_i} \partial_s D_j(r_{j,N_j}^{\nu}) & \text{if } s = \rho_{j,N_j}^{g,\nu},\, i\in\mathcal{S}_3^{\text{m}},\, j\neq i,\, \forall g,\\[13pt]
        k^{c,\nu}_{i,N_i}\, p_i^c \, \partial_s S^c_{M+1}(r_{M+1,1}^{\nu}) & \text{if } s=\rho_{M+1,1}^{g,\nu},\, i\in\mathcal{S}_2^{\text{m}},\, \forall g,\\[13pt]
        k^{c,\nu}_{i,N_i} \partial_s S_{M+1}^c(r_{M+1,1}^{\nu}) & \text{if } s=\rho_{M+1,1}^{g,\nu},\, i\in\mathcal{S}_3^{\text{m}},\, \forall g,\\[13pt]
        0 & \text{else.}
    \end{cases}
\end{equation*}
Similarly, its derivatives with respect to the control variables are
\begin{equation*}
    \partial_s \hat{\gamma}_{i,out}^{c,\nu} =
    \begin{cases}
        k^{c,\nu}_{i,N_i} \, S_{M+1}^c(r_{M+1,1}^{\nu}) & \text{if } s=p_i^c,\, i\in\mathcal{S}_2^{\text{m}},\\[5pt]
        k^{c,\nu}_{i,N_i} \partial_s D_i^c(r_{i,N_i}^{\nu}) & \text{if } s=V_{i,c},\, i\in\mathcal{S}_1^{\text{m}},\\[5pt]
        -k^{c,\nu}_{i,N_i} \partial_s D_j^c(r_{j,N_j}^{\nu}) & \text{if } s=V_{j,c},\, i\in\mathcal{S}_3^{\text{m}},\;j\neq i,\\[5pt]
        k^{c,\nu}_{i,N_i}\, p_i^c\, \partial_s S^c_{M+1}(r_{M+1,1}^{\nu}) & \text{if } s=V_{M+1,c},\, i\in\mathcal{S}_2^{\text{m}},\\[5pt]
        k^{c,\nu}_{i,N_i} \partial_s S_{M+1}^c(r_{M+1,1}^{\nu}) & \text{if } s=V_{M+1,c},\, i\in\mathcal{S}_3^{\text{m}},\\[5pt]
        0 & \text{else.}
    \end{cases}
\end{equation*}

\subsection{\texorpdfstring{$1\times M$}{1xM} FIFO diverging junction} \label{app:DAG:1xM_FIFO}

Consider a $1\times M$ diverging junction. The FIFO flux $\hat{\gamma}^{c,\nu}_{1,out}$ can be written as
\begin{equation*}
    \hat{\gamma}^{c,\nu}_{1,out} = k_{1,N_1}^{c,\nu}\, \beta, \quad \beta\coloneqq \min \left\{ D_1^c(r_{1,N_1}^{\nu}),\frac{S_2^c(r_{2,1}^{\nu})}{\alpha_2^c},\ldots,\frac{S_{M+1}^c(r_{M+1,1}^{\nu})}{\alpha_{M+1}^c} \right\}.
\end{equation*}
Denoting the index over the $M$ diverging nodes as $i=2,\ldots,M+1$, its partial derivatives with respect to the state variables would be
\begin{equation*}
    \partial_s \hat{\gamma}^{c,\nu}_{1,out} =
    \begin{cases}
        (k^{c,\nu}_{1,N_1})^{'} D_1^c(r_{1,N_1}^{\nu}) + k^{c,\nu}_{1,N_1} \partial_s D_1^c(r_{1,N_1}^{\nu}) & \text{if } s=\rho_{1,N_1}^{c,\nu},\, \beta \in\mathcal{S}^{\text{F}}_1,\\[10pt]
        k_{1,N_1}^{c,\nu}\left[-\dfrac{D_1^c(r_{1,N_1}^{\nu})}{r_{1,N_1}^{\nu}} + \partial_s D_1^c(r_{1,N_1}^{\nu})\right] & \text{if } s=\rho_{1,N_1}^{g,\nu},\, \beta \in\mathcal{S}^{\text{F}}_1,\, g\neq c,\\[10pt]
        (k^{c,\nu}_{1,N_1})^{'} \, \dfrac{S_i^c(r_{i,1}^{\nu})}{\alpha_i^c} & \text{if } s=\rho_{1,N_1}^{c,\nu},\; \beta\notin\mathcal{S}^{\text{F}}_1,\\[10pt]
        -\dfrac{k_{1,N_1}^{c,\nu}}{r_{1,N_1}^{\nu}}\cdot \dfrac{S_i^c(r_{i,1}^{\nu})}{\alpha_i^c} & \text{if } s=\rho_{1,N_1}^{g,\nu},\, \beta \notin\mathcal{S}^{\text{F}}_1,\, g\neq c\\[10pt]
        \dfrac{k^{c,\nu}_{1,N_1}}{\alpha_i^c} \, \partial_s S_i^c(r_{i,1}^{\nu}) & \text{if } s=\rho_{i,1}^{g,\nu},\; \beta \in\mathcal{S}^{\text{F}}_i,\, \forall g,\\
        0 & \text{else,}
    \end{cases}
\end{equation*}
where
\begin{align*}
    \mathcal{S}_1^{\text{F}} &\coloneqq \big\{ \beta\; :\;\beta = D_1^c(r_{1,N_1}^{\nu}) \big\},\\
    \mathcal{S}_i^{\text{F}} &\coloneqq \big\{ \beta\; :\; \beta = S_i^c(r_{i,1}^{\nu})/\alpha_i^c \big\},\\
    \mathcal{S}_i^{\text{F}}&\cap\mathcal{S}_j^{\text{F}} = \varnothing \quad\forall j\neq i,
\end{align*}
and with respect to the control variables
\begin{equation*}
    \partial_s \hat{\gamma}^{c,\nu}_{1,out} =
    \begin{cases}
        - k^{c,\nu}_{1,N_1} \dfrac{S^c_i(r_{i,1}^{\nu})}{(\alpha_i^c)^2} & \text{if } s=\alpha_i^c,\, \beta \in\mathcal{S}^{\text{F}}_i,\\[5pt]
        k^{c,\nu}_{1,N_1} \partial_s D_1^c(r_{1,N_1}^{\nu}) & \text{if } s=V_{1,c},\, \beta \in\mathcal{S}^{\text{F}}_1,\\[5pt]
        \dfrac{k^{c,\nu}_{1,N_1}}{\alpha_i^c}\, \partial_s S_i^c(r_{i,1}^{\nu}) & \text{if } s=V_{i,c},\, \beta \in\mathcal{S}^{\text{F}}_i,\\[5pt]
        0 & \text{else.}
    \end{cases}
\end{equation*}
For a $1\times M$ junction, the flux for the outgoing roads has to be differentiated as well. For the FIFO flux,
\begin{equation*}
    \hat{\gamma}_{i,in}^{c,\nu} = \alpha_i^c \hat{\gamma}_{1,out}^{c,\nu}\,, \quad i=2,\ldots,M+1.
\end{equation*}
Their partial derivatives with respect to the state variables would be, for fixed $i=2,\ldots,M+1$,

\begin{equation*}
    \partial_s \hat{\gamma}_{i,in}^{c,\nu} =
    \begin{cases}
        \alpha_i^c \left[ (k^{c,\nu}_{1,N_1})^{'} D_1^c(r_{1,N_1}^{\nu}) + k^{c,\nu}_{1,N_1} \partial_s D_1^c(r_{1,N_1}^{\nu}) \right] & \text{if } s=\rho_{1,N_1}^{c,\nu},\, \beta \in\mathcal{S}^{\text{F}}_1,\\[10pt]
        \alpha_i^c\, k_{1,N_1}^{c,\nu}\left[-\dfrac{D_1^c(r_{1,N_1}^{\nu})}{r_{1,N_1}^{\nu}} + \partial_s D_1^c(r_{1,N_1}^{\nu})\right] & \text{if } s=\rho_{1,N_1}^{g,\nu},\, \beta \in\mathcal{S}^{\text{F}}_1,\, g\neq c,\\[10pt]
        \dfrac{\alpha_i^c}{\alpha_j^c} \, (k^{c,\nu}_{1,N_1})^{'} S_j^c(r_{j,1}^{\nu}) & \text{if } s=\rho_{1,N_1}^{c,\nu},\, \beta \in\mathcal{S}^{\text{F}}_j,\, \forall j,\\[8pt]
        -\dfrac{\alpha_i^c}{\alpha_j^c} \cdot \dfrac{k^{c,\nu}_{1,N_1}}{r_{1,N_1}^{\nu}} \, S_j^c(r_{j,1}^{\nu}) & \text{if } s=\rho_{1,N_1}^{g,\nu}, \, \beta \in\mathcal{S}^{\text{F}}_j,\, \forall j,\, g\neq c,\\[10pt]
        \dfrac{\alpha_i^c}{\alpha_j^c} \, k^{c,\nu}_{1,N_1} \partial_s S_i^c(r_{i,1}^{\nu}) & \text{if } s=\rho_{j,1}^{g,\nu},\, \beta \in\mathcal{S}^{\text{F}}_j,\, \forall j,\, \forall g,\\
        0 & \text{else,}
    \end{cases}
\end{equation*}
and with respect to the control variables
\begin{equation*}
    \partial_s \hat{\gamma}_{i,in}^{c,\nu} =
    \begin{cases}
        k^{c,\nu}_{1,N_1} D_1^c(r_{1,N_1}^{\nu}) & \text{if } s=\alpha_i^c,\, \beta \in\mathcal{S}^{\text{F}}_1,\\[5pt]
        k^{c,\nu}_{1,N_1} \dfrac{S_j^c(r_{j,1}^{\nu})}{\alpha_j^c} & \text{if } s=\alpha_i^c,\, \beta \in\mathcal{S}^{\text{F}}_j,\, j\neq i,\\[5pt]
        \alpha_i^c\, k^{c,\nu}_{1,N_1} \partial_s D_1^c(r_{1,N_1}^{\nu}) & \text{if } s=V_{1,c},\, \beta \in\mathcal{S}^{\text{F}}_1,\\[5pt]
        \dfrac{\alpha_i^c}{\alpha_j^c}\, k^{c,\nu}_{1,N_1} \partial_s S_j^c(r_{j,1}^{\nu}) & \text{if } s=V_{j,c},\, \beta \in\mathcal{S}^{\text{F}}_j,\, \forall j,\\
        0 & \text{else.}
    \end{cases}
\end{equation*}

\subsection{\texorpdfstring{$1\times M$}{1xM} non-FIFO diverging junction} \label{app:DAG:1xM_nonFIFO}

The non-FIFO flux for link $\hat{\gamma}_{1,out}^{c,\nu}$ of a $1\times M$ junction can be written as the sum of the links $\hat{\gamma}_{i,in}^{c,\nu}$, $i=2,\ldots,M+1$,
\begin{equation*}
    \hat{\gamma}_{1,out}^{c,\nu} = \sum_{i=2}^{M+1} \hat{\gamma}_{i,in}^{c,\nu}, \quad \hat{\gamma}_{i,in}^{c,\nu} = k_{1,N_1}^{c,\nu} \min \left\lbrace \alpha_i^c D_1^c(r_{1,N_1}^{\nu}),\, S_i^c(r_{i,1}^{\nu}) \right\rbrace.
\end{equation*}
Its partial derivatives with respect to $s$ would therefore be
\begin{equation*}
    \partial_s \hat{\gamma}^{c,\nu}_{1,out} = \sum_{i=2}^{M+1}\partial_s \hat{\gamma}^{c,\nu}_{i,in}.
\end{equation*}
The derivatives with respect to the state variables of the outgoing non-FIFO fluxes are, for $i=2,\ldots,M+1$,

\begin{equation*}
    \partial_s \hat{\gamma}^{c,\nu}_{i,in} =
    \begin{cases}
        \alpha_i^c \left[ (k^{c,\nu}_{1,N_1})^{'} D_1^c(r_{1,N_1}^{\nu}) + k^{c,\nu}_{1,N_1} \partial_s D_1^c(r_{1,N_1}^{\nu}) \right] & \text{if } s=\rho_{1,N_1}^{c,\nu},\, i\in\mathcal{S}_1^{\text{nF}},\\[9pt]
        \alpha_i^c \, k^{c,\nu}_{1,N_1}\left[ -\dfrac{D_1^c(r_{1,N_1}^{\nu})}{r_{1,N_1}^{\nu}} + \partial_s D_1^c(r_{1,N_1}^{\nu}) \right] & \text{if } s=\rho_{1,N_1}^{g,\nu},\, i\in\mathcal{S}_1^{\text{nF}},\, g\neq c,\\[13pt]
        (k^{c,\nu}_{1,N_1})^{'} S_i^c(r_{i,1}^{\nu}) & \text{if } s=\rho_{1,N_1}^{c,\nu},\, i\in\mathcal{S}_2^{\text{nF}},\\[9pt]
        -\dfrac{k^{c,\nu}_{1,N_1}}{r_{1,N_1}^{\nu}} \, S_i^c(r_{i,1}^{\nu}) & \text{if } s=\rho_{1,N_1}^{g,\nu},\, i\in\mathcal{S}_2^{\text{nF}},\, g\neq c,\\[11pt]
        k^{c,\nu}_{1,N_1} \partial_s S_i^c(r_{i,1}^{\nu}) & \text{if } s=\rho_{i,1}^{g,\nu},\, i\in\mathcal{S}_2^{\text{nF}},\, \forall g,\\[9pt]
        0 & \text{else,}
    \end{cases}
\end{equation*}
where
\begin{align*}
    \mathcal{S}_1^{\text{nF}} &\coloneqq \big\{ i\; :\; \alpha_i^c D_1^c(r_{1,N_1}^{\nu})\leq S_i^c(r_{i,1}^{\nu}) \big\},\\
    \mathcal{S}_2^{\text{nF}} &\coloneqq \big\{ i\; :\; \alpha_i^c D_1^c(r_{1,N_1}^{\nu})> S_i^c(r_{i,1}^{\nu}) \big\},\\
    \mathcal{S}_1^{\text{nF}}&\cap\mathcal{S}_2^{\text{nF}} = \varnothing,
\end{align*}
and with respect to the control variables,
\begin{equation*}
    \partial_s \hat{\gamma}^{c,\nu}_{i,in} =
    \begin{cases}
        k^{c,\nu}_{1,N_1} D_1^c(r_{1,N_1}^{\nu}) & \text{if } s=\alpha_i^c,\, i\in\mathcal{S}_1^{\text{nF}},\\[5pt]
        k^{c,\nu}_{1,N_1} \alpha_i^c \, \partial_s D_1^c(r_{1,N_1}^{\nu}) & \text{if } s=V_{1,c},\, i\in\mathcal{S}_1^{\text{nF}},\\[5pt]
        k^{c,\nu}_{1,N_1} \partial_s S_i^c(r_{i,1}^{\nu}) & \text{if } s=V_{i,c},\, i\in\mathcal{S}_2^{\text{nF}},\\[5pt]
        0 & \text{else.}
    \end{cases}
\end{equation*}

\subsection{Boundary conditions} \label{app:DAG:BC}

Focusing on the derivatives for the boundary conditions, the inflows $\hat{\gamma}_{\ell,in}^{c,\nu}$, $\ell\in\mathcal{O}$, are defined as
\begin{align*}
    \hat{\gamma}_{\ell,in}^{c,\nu} = \min \bigg\{ D_{\ell}^c(l_{\ell}^{c,\nu}),\, \max \big\{ \frac{1}{N} S_{\ell}^c(r_{\ell,1}^{\nu}),\, S_{\ell}^c(r_{\ell,1}^{\nu}) - \sum_{g\neq c} D_{\ell}^g(l_{\ell}^{g,\nu}) \big\} \bigg\}.
\end{align*}
Defining
\begin{align*}
    \mathcal{S}_1^{in} &\coloneqq \big\{ \ell\; :\; \ell\in\mathcal{O},\; \hat{\gamma}_{\ell,in}^{c,\nu} = D_{\ell}^c(l_{\ell}^{c,\nu}) \big\},\\[3pt]
    \mathcal{S}_2^{in} &\coloneqq \big\{ \ell\; :\; \ell\in\mathcal{O},\; \hat{\gamma}_{\ell,in}^{c,\nu} = S_{\ell}^c(r_{\ell,1}^{\nu})/N \big\},\\[3pt]
    \mathcal{S}_3^{in} &\coloneqq \big\{ \ell\; :\; \ell\in\mathcal{O},\; \hat{\gamma}_{\ell,in}^{c,\nu} = S^c_{\ell}(r_{\ell,1}^{\nu}) - \sum_{g\neq c} D^g_j(l_{\ell}^{g,\nu}) \big\},\\
    \mathcal{S}_i^{in}&\cap\mathcal{S}_j^{in} = \varnothing \quad\forall j\neq i,
\end{align*}
the partial derivatives with respect to the state variables would therefore be
\begin{equation*}
    \partial_s \hat{\gamma}_{\ell,in}^{c,\nu} =
    \begin{cases}
        \dfrac{1}{N}\,\partial_s S_{\ell}^c(r_{\ell,1}^{\nu}) & \text{if } s=\rho_{\ell,1}^{g,\nu},\, \ell\in\mathcal{S}^{in}_2,\, \forall g,\\[5pt]
        \partial_s S_{\ell}^c(r_{\ell,1}^{\nu}) & \text{if } s=\rho_{\ell,1}^{g,\nu},\, \ell\in\mathcal{S}^{in}_3,\, \forall g,\\[5pt]
        0 & \text{ else,}
    \end{cases}
\end{equation*}
and with respect to the control variables
\begin{equation*}
    \partial_s \hat{\gamma}_{\ell,in}^{c,\nu} =
    \begin{cases}
        \partial_s D_{\ell}^c(l_{\ell}^{c,\nu}) & \text{if } s=V_{\ell,c},\, \ell\in\mathcal{S}^{in}_1,\\[6pt]
        \dfrac{1}{N}\, \partial_s S_{\ell}^c(r_{\ell,1}^{\nu}) & \text{if } s=V_{\ell,c},\, \ell\in\mathcal{S}^{in}_2,\\[6pt]
        \partial_s S_{\ell}^c(r_{\ell,1}^{\nu}) & \text{if } s=V_{\ell,c},\, \ell\in\mathcal{S}^{in}_3,\\
        0 & \text{else,}
    \end{cases}
\end{equation*}
where
\begin{align*}
    \partial_{V_{\ell,c}} D_{\ell}^c(l_{\ell}^{c,\nu}) &=
    \begin{cases}
        \partial_{V_{\ell,c}} Q_{\ell}^c(r_{\ell,cr}^c) & \text{if } l_{\ell}^{c,\nu}>0,\\
        \partial_{V_{\ell,c}} F_{\ell,in}^c & \text{if } l_{\ell}^{c,\nu}=0,
    \end{cases}\\[3pt]
    &=
    \begin{cases}
        r_{\ell,cr}^c \partial_{V_{\ell,c}} v_{\ell,c}(r_{\ell,cr}^c) & \text{if } l_{\ell}^{c,\nu}>0,\\
        0 & \text{if } l_{\ell}^{c,\nu}=0.
    \end{cases}
\end{align*}
For the outflows $\hat{\gamma}_{\ell,out}^{c,\nu}$ where $\ell$ is the index corresponding to a road leaving the network,
\begin{equation*}
    \hat{\gamma}_{\ell,out}^{c,\nu} = \min\left\{ k^{c,\nu}_{\ell,N_{\ell}} D_{\ell}^c(r_{\ell,N_{\ell}}^{\nu}),\, F^c_{\ell,out} \right\}.
\end{equation*}
Their partial derivatives with respect to the state variables are
\begin{equation*}
    \partial_s \hat{\gamma}_{\ell,out}^{c,\nu} =
    \begin{cases}
        (k^{c,\nu}_{\ell,N_{\ell}})^{'} D_{\ell}^c(r_{\ell,N_{\ell}}^{\nu}) + k^{c,\nu}_{\ell,N_{\ell}} \partial_s D_{\ell}^c(r_{\ell,N_{\ell}}^{\nu}) & \text{if } s=\rho_{\ell,N_{\ell}}^{c,\nu},\, \ell\in\mathcal{S}^{out},\\[7pt]
        k_{\ell,N_{\ell}}^{c,\nu}\left[-\dfrac{D_{\ell}^c(r_{\ell,N_{\ell}}^{\nu})}{r_{\ell,N_{\ell}}^{\nu}} + \partial_s D_{\ell}^c(r_{\ell,N_{\ell}}^{\nu})\right] & \text{if } s=\rho_{\ell,N_{\ell}}^{g,\nu},\, \ell\in\mathcal{S}^{out},\, g\neq c,\\[9pt]
        0 & \text{else,}
    \end{cases}
\end{equation*}
where
\begin{equation*}
    \mathcal{S}^{out} \coloneqq \big\{ \ell\; :\; k^{c,\nu}_{\ell,N_{\ell}}D_{\ell}^c(r_{\ell,N_{\ell}}^{\nu})\leq F_{\ell,out}^c \big\},
\end{equation*}
while with respect to the control variables they are
\begin{equation*}
    \partial_s \hat{\gamma}_{\ell,out}^{c,\nu} =
    \begin{cases}
        k^{c,\nu}_{\ell,N_{\ell}} \partial_s D_{\ell}^c(r_{\ell,N_{\ell}}^{\nu}) & \text{if } s=V_{\ell,c},\, \ell\in\mathcal{S}^{out},\\
        0 & \text{else.}
    \end{cases}
\end{equation*}

\end{appendices}


\bibliography{bibarticleDAG.bib}% common bib file

\end{document}